\documentclass[preprint,number,sort&compress,12pt]{elsarticle}

\usepackage[hyphens]{url}
\usepackage{graphicx}
\usepackage{hyperref}
\usepackage{amsmath,amssymb,amsfonts}
\usepackage{bm}
\usepackage{eucal}
\usepackage{color}
\usepackage{booktabs}
\usepackage{tabularx}
\usepackage{etoolbox}
\usepackage{tikz}
\usepackage[margin=2.5cm]{geometry}

\makeatletter
\apptocmd{\ps@pprintTitle}{%
  \apptocmd{\@oddfoot}{\unskip}{}{}%
  \let\@evenfoot\@oddfoot
}{}{}
\makeatother
\hypersetup{
	colorlinks   = true, 
	urlcolor     = blue, 
	linkcolor    = blue, 
	citecolor   = red 
}
\usetikzlibrary{calc}
\usetikzlibrary{shapes}
\usetikzlibrary{arrows.meta,positioning,fit,backgrounds}

\usepackage{fancyvrb}
\usepackage{tcolorbox}
\tcbuselibrary{skins,breakable}
\input{listings/pygments_style}
\newtcolorbox{codebox}[1][]{
  breakable, enhanced, sharp corners=downhill,
  colback=gray!4, colframe=gray!45, boxrule=0.4pt,
  left=3pt, right=3pt, top=2pt, bottom=2pt,
  fontupper=\footnotesize, #1}

\begin{document}

\begin{frontmatter}
		
		\title{JAXIGA - A software framework for GPU-enabled and differentiable isogeometric analysis}

		\author[cmu]{Cosmin Anitescu \fnref{fn1}}
        \author[cmu]{Yongjie Jessica Zhang \corref{cor1}\fnref{fn2}}
		
		\cortext[cor1]{Corresponding author}
		\address[cmu]{Department of Mechanical Engineering, Carnegie Mellon University, Pittsburgh, PA, USA}
		
		\fntext[fn1]{canitesc@andrew.cmu.edu}
		\fntext[fn2]{jessicaz@andrew.cmu.edu}

\begin{abstract}
In this paper, we present JAXIGA, an open-source library for isogeometric
analysis based on JAX. The library separates the construction of the mesh from
the computations which are traced for compilation and differentiation. This
allows geometry evaluation, assembly and the supported solvers to be combined in
a just-in-time compiled program, for which first-order reverse-mode
derivatives can be computed. We consider three numerical methods obtained
from a single pointwise energy density. The Galerkin and energy-minimization
methods use the same discrete potential and have the same stationary points,
while collocation enforces the corresponding strong form. The sensitivities
of a converged solution are obtained by solving one adjoint problem, and their
accuracy depends on the tolerances of both the primal and adjoint problems.
We also describe hierarchical extraction for local refinement of the
variational methods and the use of smooth spline spaces for fourth-order
problems. The numerical results are compared with analytical solutions and
independent FEniCSx and JAX-FEM computations. For the reported GPU benchmark with
$36{,}448$ unknowns, assembly is twenty-four times faster than on the host
CPU. When the conjugate-gradient solver is also run on the GPU, the
total solution time is reduced from $31.8$ to $1.38$ seconds, and the result
agrees with the direct solution to $5\times10^{-11}$. These results show that
the choice of solver and the reuse of compilation are important for obtaining
an overall reduction in computational cost.

\end{abstract}


\end{frontmatter}

\section{Introduction}
\label{sec:intro}

Isogeometric Analysis (IGA) uses spline basis functions to represent the
geometry and to approximate the solution of a partial differential equation
(PDE)~\cite{hughes2005,cottrell2009}. One of the advantages of this approach is
that NURBS can represent conic geometry exactly. Moreover, with suitable
degrees and knot multiplicities, the basis functions have the $C^1$
continuity required by weak formulations of fourth-order PDEs and by
collocation formulations of second-order PDEs. Higher continuity can also
reduce the number of degrees of freedom as the polynomial degree increases.
These properties are useful in inverse problems and design optimization,
provided that sensitivities can be computed through the evaluation of the
geometry, the assembly procedure and the solution of the discrete equations.

In JAXIGA, the PDE and boundary conditions are expressed using
JAX~\cite{jax2018}. For a fixed discretization, the user specifies a pointwise
energy density, from which the Galerkin, energy-minimization and collocation
methods are obtained. First-order reverse-mode sensitivities of the converged
solutions can then be computed, with GPU acceleration used where available.
We treat mesh construction and refinement separately from these computations.
As will be discussed below, this separation determines both the capabilities
of the library and the conditions under which differentiation is supported.

\subsection{Related work and software}
\label{subsec:existing_software}

Several libraries have been developed for isogeometric analysis, with different
choices of geometry representation, programming environment and solver.
PetIGA~\cite{petiga}
provides PETSc-based parallel analysis, GeoPDEs~\cite{geopdes,geopdeshier} supports
MATLAB/Octave prototyping and hierarchical spline adaptivity, and
G+Smo~\cite{gismo} provides extensive spline geometry and analysis facilities. Among Python-based packages, tIGAr~\cite{tigar} connects spline extraction to FEniCS, while
Nutils~\cite{nutils} offers a finite element environment with spline
support. FEniCSx~\cite{dolfinx}, Firedrake~\cite{firedrake} and
deal.II~\cite{dealii} provide broader discretization and solver infrastructure.
Sensitivities can be obtained through a separate adjoint framework, such as
\texttt{pyadjoint}~\cite{pyadjoint}.

The design of JAXIGA is closely related to that of JAX-FEM~\cite{jaxfem}, which
uses local constitutive functions, structured parameter data and implicit
differentiation in JAX. We follow some of its interface conventions, while
using NURBS and hierarchical splines in place of Lagrange elements.
Differentiable fluid and dynamics frameworks
\cite{jaxcfd,jaxfluids,phiflow,diffrax,hu2019} demonstrate the wider applicability
of this programming model. The energy method is also related to the deep energy
method~\cite{samaniego2020}, but uses a spline trial space rather than a neural
network.

The use of GPUs to accelerate IGA assembly has also been studied previously.
Karatarakis et al.~\cite{karatarakis2014} developed an interaction-wise GPU
procedure for stiffness assembly and reported orders-of-magnitude reductions
in assembly time relative to standard element-wise Gauss quadrature in their
benchmarks. Kiran et al.~\cite{kiran2019} implemented GPU-based FE and IGA
computations within MATLAB. Their maximum $41.4\times$ speedup over a
single-core sequential CPU implementation concerned numerical integration and
assembly for a linear FE benchmark. They also demonstrated nonlinear IGA for
an adhesive strip-peeling problem. We note that these measurements concern
different problems, implementations and hardware from those considered here.
Therefore, they do not provide a direct comparison of the total solution cost
with JAXIGA.

Efficient IGA also depends on the assembly algorithm, independently of the
programming framework. Established approaches include weighted
quadrature~\cite{calabro2017}, sum factorization for assembled and matrix-free
operators~\cite{bressan2019}, and low-rank tensor assembly~\cite{lowrank2015}.
The Python toolbox pyiga~\cite{pyiga} provides generated assembly kernels,
low-rank assembly and hierarchical spline adaptivity. Sangalli and
Tani~\cite{sangalli2018} combined weighted quadrature with matrix-free operator
application and preconditioning to reduce both runtime and storage. These
methods provide possible alternatives to the quadrature and assembly
procedures currently used in JAXIGA, and could be considered in further
development of the library.

Automatic differentiation has also been used to compute sensitivities in
spline-based design.
L\'opez, Anitescu and Rabczuk~\cite{lopez2021} used forward- and reverse-mode
automatic differentiation for multilevel isogeometric shape optimization,
propagating sensitivities between refined analysis meshes and a coarse design
parameterization. Wu~\cite{wu2023} combined JAX automatic differentiation,
adjoint sensitivities and XLA just-in-time compilation for structural shape
optimization, including B\'ezier-surface parameterizations. Thus, the use of
automatic differentiation for IGA shape sensitivities, and the use of JAX for
adjoint shape optimization, have precedents in the literature.

A further related framework is IgANet~\cite{iganet,iganets2026}, a C++ library
built on LibTorch that combines B-splines and NURBS with PDE solving and deep
operator learning. Its implementation supports automatic differentiation of
spline data, GPU execution and runtime compilation of spline specializations.
These capabilities overlap with JAXIGA's use of differentiable spline
representations on accelerators, although the IGANets approach focuses on
learning maps between problem data and spline solution coefficients.

In Table~\ref{tab:comparison}, we compare the capabilities available within
several representative frameworks. Some of these provide more extensive mesh
handling, element libraries, adaptivity, solvers and parallel scaling than
JAXIGA. Here we consider a smaller set of problem classes, for which the
built-in methods can be extended through a common interface.

\begin{table}[htbp]
\centering
\caption{Capability comparison. ``Autodiff through solver'' means gradients of
the converged solution are available in the framework's own language, rather
than through a separate adjoint layer. ``Native GPU'' means the same source runs
on an accelerator without a separate implementation. For multipatch coupling,
``automatic'' denotes built-in construction of interface continuity constraints
for supported conforming patches, with interface topology supplied or detected.
``Manual'' requires user-supplied coupling code or constraints. HB and THB
denote hierarchical and truncated hierarchical B-splines. ``Unverified'' means
that the capability was not established from the cited sources, rather than
that it is absent.}
\label{tab:comparison}
\small
\begin{tabular}{l l c c c c}
\hline
Framework & Language & Discretization & Autodiff & Native & Multipatch \\
 & & & through solver & GPU & coupling \\
\hline
PetIGA~\cite{petiga}   & C       & NURBS    & no      & partial & manual \\
GeoPDEs~\cite{geopdes,geopdeshier} & MATLAB  & NURBS/HB/THB & no & no & automatic \\
G+Smo~\cite{gismo}     & C++     & NURBS/THB& no      & no      & automatic \\
tIGAr~\cite{tigar}     & Python/C++ & NURBS & via adjoint layer & no & manual \\
FEniCSx~\cite{dolfinx}  & Python/C++ & FEM   & via \texttt{pyadjoint} & partial & n/a \\
Firedrake~\cite{firedrake} & Python & FEM   & via \texttt{pyadjoint} & partial & n/a \\
JAX-FEM~\cite{jaxfem}  & Python  & FEM      & yes     & yes     & n/a \\
IgANet~\cite{iganet,iganets2026} & C++ & NURBS & unverified & yes & unverified \\
\textbf{JAXIGA}        & Python  & NURBS/THB & yes     & yes     & automatic \\
\hline
\end{tabular}
\end{table}

\subsection{Contribution and scope}
\label{subsec:contributions}

The main contribution of this work is a common JAX interface in which three
numerical methods are obtained from one energy definition, together with
implicit solution sensitivities, conforming multipatch coupling and
hierarchical spline extraction. We keep the static connectivity separate from
the geometry and field data which are traced by JAX. In particular, this
allows evaluation, assembly and the supported solvers to be combined under
\texttt{jax.jit} and differentiated with \texttt{jax.grad}, provided that the
discrete solution is converged and regular. The resulting sensitivities are
first-order reverse-mode derivatives. We also use the same extraction
interface for multipatch coupling and truncated hierarchical B-splines.
Therefore, the variational methods can be applied to locally refined spaces
without a separate implementation of the computation. Collocation, however,
is restricted to tensor-product spaces in the current version.

In the numerical experiments, we examine the accuracy of the discretization
and the computational cost of the implementation. The computed solutions are
compared with analytical solutions and independent FEniCSx and JAX-FEM
results. We also compare the sensitivities with finite-difference
approximations in order to verify the differentiation procedure.

The remainder of the paper is organized as follows.
Section~\ref{sec:math} describes the spline discretization and the three
numerical methods. In Section~\ref{sec:design}, we present the software design
and the conditions required for differentiation. Numerical results for the
accuracy, sensitivities and computational cost are given in
Section~\ref{sec:examples}. Finally, the limitations of the current
implementation are discussed in Section~\ref{sec:discussion}.

\section{Mathematical Preliminaries}
\label{sec:math}

In this section, we first describe the spline representation of the geometry
and the approximate solution. We then use this representation to discretize a
variational problem, from which the three numerical methods and their
solution strategies will be obtained.

\subsection{Spline discretization}
\label{subsec:nurbs}

Let $N_{i,p}$ be the degree-$p$ B-spline basis on an open knot vector
$\Xi=(\xi_1,\ldots,\xi_{n+p+1})$, defined by the Cox--de Boor
recursion~\cite{piegl1997}. The basis is nonnegative, forms a partition of
unity on $[\xi_{p+1},\xi_{n+1}]$, and has $C^{p-k}$ continuity at an interior
knot of multiplicity $k\leq p$.

In several parametric dimensions, tensor products of the univariate functions
give the basis functions $B_A$, where $A$ indexes the tensor-product basis.
Assigning a positive weight $w_A$ to each function gives the rational basis,
geometry map and scalar solution approximation
\begin{equation}
 R_A(\boldsymbol\xi)=\frac{w_A B_A(\boldsymbol\xi)}
 {\sum_B w_B B_B(\boldsymbol\xi)},\qquad
 \bm{x}(\boldsymbol\xi)=\sum_A R_A(\boldsymbol\xi)\bm{P}_A,
 \qquad u_h(\bm{x})=\sum_A R_A(\boldsymbol\xi(\bm{x}))c_A .
\label{eq:nurbs}
\end{equation}
Here, $\bm P_A$, $w_A>0$ and $c_A$ are the control points, weights and
solution coefficients, and $\boldsymbol\xi(\bm x)$ is the inverse geometry map.
For vector-valued fields, the same expansion is used for each component. In
the following, we assume that the geometry map is regular, with a nonsingular
Jacobian.

Both knot insertion and degree elevation preserve this geometry map. In
particular, elevating the degree before inserting new knots gives the usual
maximally smooth $k$-refined spaces. For elementwise evaluation,
B\'ezier extraction~\cite{borden2011} expresses each active spline basis
function in an element-local Bernstein basis. This representation also
accommodates truncated hierarchical B-spline (THB) refinement, as described in
Section~\ref{subsec:design_hierarchical}.

\subsection{Variational formulation}
\label{subsec:energy_framework}
\label{subsec:energy_pdEs}

For a conservative problem, we write the total potential as
\begin{equation}
 \mathcal E(\bm u;\bm\theta)=\int_\Omega
 \psi(\nabla\bm u,\bm u,\bm{x},\bm\theta)\,d\Omega
 -\int_\Omega \bm f\cdot\bm u\,d\Omega
 -\int_{\Gamma_N}\bm t\cdot\bm u\,d\Gamma ,
\label{eq:energy_decomposition}
\end{equation}
Here, $\psi$ is the internal energy density, $\bm\theta$ collects model
parameters, and $\bm f$ and $\bm t$ are prescribed volume loads and Neumann
tractions. Derivatives are taken in physical coordinates. For scalar fields,
the load terms use ordinary multiplication.

To impose the Dirichlet data, we use a lifting, which is a field satisfying
the prescribed boundary values. This field is added to the unknown field with
homogeneous Dirichlet data.
Substituting the spline approximation into~\eqref{eq:energy_decomposition}
and evaluating the integrals by Gauss quadrature then gives a discrete potential
$\Pi(\bm d,\bm\theta)$, where $\bm d$ collects the free solution coefficients
after the Dirichlet constraints have been imposed. Differentiating this
potential gives
\begin{equation}
 \bm r=\nabla_{\bm d}\Pi,\qquad
 \bm A=\nabla^2_{\bm d}\Pi ,
\label{eq:energy_hessian}
\end{equation}
where $\bm r$ is the discrete residual vector and $\bm A$ is the tangent
matrix. The operators $\nabla_{\bm d}$ and $\nabla^2_{\bm d}$ denote the gradient and Hessian
with respect to the free coefficients. For a quadratic potential, this reduces
to $\bm r=\bm K\bm d-\bm f_h$, where $\bm K$ is the stiffness matrix and
$\bm f_h$ is the discrete load vector, including contributions from prescribed
Dirichlet values. Consequently, finding a stationary point requires only the
solution of a linear system. For a nonconvex potential, however, stationarity does not imply
minimality, and Newton iteration and energy minimization may reach different
stationary points.

Two representative choices of energy density are those for Poisson's equation
and small-strain linear elasticity:
\begin{align}
 \psi_{\mathrm{Poisson}}&=\tfrac12 a(\bm x)|\nabla u|^2,
 &-\nabla\cdot(a\nabla u)&=f,\quad a>0,\label{eq:poisson_energy}\\
 \psi_{\mathrm{elasticity}}&=\tfrac12\bm\varepsilon:\mathbb C:\bm\varepsilon,
 &\bm\varepsilon&=\tfrac12(\nabla\bm u+\nabla\bm u^{\!\top}).
 \label{eq:elasticity_energy}
\end{align}
The coefficient $a$ is positive, $\bm\varepsilon$ is the infinitesimal strain,
$\mathbb C$ is the elasticity tensor, and the colon denotes tensor contraction.

The same formulation can be used for nonlinear material models by replacing
$\psi$ with the corresponding energy density. For fourth-order variational
models, the energy also depends
on $\nabla^2u$, the Hessian of the field in physical coordinates. These models
require an $H^2$-conforming trial space, whose functions and weak derivatives through
order two are square-integrable. The history-dependent fracture formulation
requires additional state variables and is therefore specified separately in
Section~\ref{subsec:ex_phase_field}.

\subsection{Methods and coupling}
\label{subsec:collocation}
\label{subsec:multipatch}

We now consider the three numerical methods. The Galerkin method seeks a
stationary point by solving $\bm r=\bm0$, while energy minimization searches
directly for a minimum of $\Pi$. Thus, the two methods use the same discrete
potential. For a scalar problem,
collocation~\cite{auricchio2010,schillinger2013} enforces the associated strong
form,
$-\nabla\cdot(\partial\psi/\partial\nabla u)+\partial\psi/\partial u=f$
at Greville points. Here, $\partial\psi/\partial\nabla u$ is the flux derived
from the energy density, and $\partial\psi/\partial u$ is its derivative with
respect to the field.
For each univariate basis function, the Greville coordinate is the knot average
$(\xi_{i+1}+\cdots+\xi_{i+p})/p$. Taking tensor products of these coordinates
and applying the geometry map gives the physical collocation points. This
construction avoids volume quadrature in the collocation equations. However,
the cost and convergence of the method also depend on the arrangement of the
equations and on the linear solver. Adding
Neumann and interface flux equations generally makes the system overdetermined,
so the implementation solves it in the least-squares sense.
Table~\ref{tab:method_scope} summarizes the resulting discrete problems, and
Section~\ref{subsec:one_definition} describes their implementation.

\begin{table}[htbp]
\centering
\caption{Methods derived from the same energy density. $\bm F$ is the vector
of strong-form, Neumann and interface residuals, and $\bm J=\partial\bm F/\partial\bm d$.
The collocation stationarity condition generally permits $\bm F\ne\bm0$.}
\label{tab:method_scope}
\small
\begin{tabularx}{\textwidth}{l X X}
\toprule
Method & Discrete problem & Supported spaces \\
\midrule
Galerkin & $\nabla_{\bm d}\Pi=\bm0$ & Tensor-product and THB \\
Energy minimization & $\min_{\bm d}\Pi$ & Tensor-product and THB \\
Collocation & $\min_{\bm d}\tfrac12\|\bm F\|^2$;
$\bm J^{\!\top}\bm F=\bm0$ & Tensor-product only \\
\bottomrule
\end{tabularx}
\end{table}

To extend the discretization to several patches, conforming coupling identifies
matching degrees of freedom in their traces, that is, the restrictions of the
fields to the interface. This imposes $u_\alpha=u_\beta$ on the common boundary,
where $\alpha$ and $\beta$ label neighboring patches and $u_\alpha$ and
$u_\beta$ are their respective solution fields. The total potential is then
the sum of the patch contributions. Such matching requires compatible trace
bases and orientations, in addition to coincident control points. These
requirements are detailed in Section~\ref{subsec:design_multipatch}. The resulting
coupling provides $C^0$ continuity, meaning continuity of the field across the
interface. Nonconforming coupling and $C^1$ coupling, which would also enforce
continuity of first derivatives, are not implemented in the current version.

\subsection{Solution strategies}
\label{subsec:solution_strategies}

Linear Galerkin systems are solved by a sparse LU factorization (SuperLU) on
the host CPU, a dense direct solver for small systems, or a Krylov method with optional
Jacobi (diagonal) preconditioning. Automatic selection of the
conjugate-gradient (CG) method requires a linear problem that declares a
positive-definite operator, admissible material data and constraints that
remove its nullspace. Within these restrictions the choice depends on the
backend and system size. On the GPU, CG is selected above about $3{,}000$ free
unknowns, while SuperLU is used on the CPU and for smaller systems. This
threshold was calibrated on the plate benchmark of
Section~\ref{subsec:ex_parity} and may not hold for other degrees or hardware. Users can also select a solver
explicitly.

Nonlinear Galerkin problems use Newton iteration with a backtracking line
search. Energy minimization uses the \texttt{optax} implementation of the
limited-memory Broyden--Fletcher--Goldfarb--Shanno (L-BFGS)
method~\cite{liu1989} with a zoom line search, optionally initialized by a few
Adam steps, and collocation uses Gauss--Newton iteration for its least-squares
problem. Stopping conditions and implicit sensitivities are described in
Section~\ref{subsec:design_solvers}.

\section{Software Design}
\label{sec:design}

\begin{figure}[htbp]
\centering
\begin{tikzpicture}[
  font=\small,
  mod/.style   = {rectangle, rounded corners=2pt, draw, thick, align=center,
                  minimum height=8mm, inner sep=4pt},
  setupbox/.style = {mod, fill=blue!8,   draw=blue!55},
  tracedbox/.style= {mod, fill=orange!12, draw=orange!70},
  bothbox/.style  = {mod, fill=green!10,  draw=green!55!black},
  lbl/.style   = {font=\scriptsize\itshape, text=black!55},
  flow/.style  = {-{Latex[length=2mm]}, thick, black!65},
  back/.style  = {-{Latex[length=2mm]}, thick, red!65!black, dashed},
  loop/.style  = {-{Latex[length=2mm]}, thick, blue!60!black, densely dotted},
]

\useasboundingbox (-0.45, -9.35) rectangle (15.3, 2.25);
\fill[blue!5,   rounded corners] (-0.3,  0.62) rectangle (13.9, 1.02);
\fill[orange!6, rounded corners] (-0.3, -3.05) rectangle (13.9, -2.65);
\node[lbl, anchor=west] at (-0.15,  0.82) {setup stage --- NumPy, host, runs once, static};
\node[lbl, anchor=west] at (-0.4, -2.85)
      {compute stage: traced JAX; one density $\psi$, three methods};

\node[setupbox] (geom)  at (1.25, 0)  {\texttt{geometry}\\[-1pt]{\scriptsize Patch, primitives}};
\node[setupbox] (multi) at (4.25, 0)  {\texttt{multipatch}\\[-1pt]{\scriptsize interface detection}};
\node[setupbox] (space) at (7.4,  0)  {\texttt{space}\\[-1pt]{\scriptsize extraction, IEN, dofs}};
\node[setupbox] (pts)   at (11.95, 0)  {\texttt{pointset}\\[-1pt]{\scriptsize Gauss, Greville, grid}};
\node[setupbox] (hier)  at (11.95, -1.6)
  {\texttt{space.hierarchical}\\[-1pt]{\scriptsize THB local refinement}};
\draw[flow] (pts.south) -- (hier.north);

\draw[flow] (geom)  -- (multi);
\draw[flow] (multi) -- (space);
\draw[flow] (space) -- (pts);

\node[tracedbox] (prob) at (1.9, -1.6)
  {\texttt{forms}\\[-1pt]{\scriptsize $\psi(\nabla u, u, \bm{x}, \bm{\theta})$, BCs}};
\node[bothbox]  (eval) at (7.4, -1.6)
  {\texttt{space.evaluation}\\[-1pt]{\scriptsize $R,\ \nabla R,\ \nabla^2 R,\ \bm{x},\ w,\ \bm{n}$}};

\draw[flow] (space.south) -- (eval.north);
\draw[flow] (pts.south) -- (eval.north east);
\draw[flow] (hier.west) -- (eval.east);
\draw[flow] (prob.east) -- (eval.west);

\node[tracedbox] (gal) at (2.0, -4.0)
  {\texttt{methods.galerkin}\\[-1pt]{\scriptsize $\partial\Pi/\partial\bm{d}=\bm{0}$}};
\node[tracedbox] (ene) at (6.1, -4.0)
  {\texttt{methods.energy}\\[-1pt]{\scriptsize $\min\ \Pi$}};
\node[tracedbox] (col) at (10.2, -4.0)
  {\texttt{methods.collocation}\\[-1pt]{\scriptsize strong form, Greville}};

\draw[thick, black!65] (eval.south) -- (7.4,-3.2);
\draw[thick, black!65] (2.0,-3.2) -- (10.2,-3.2);
\draw[flow] (2.0,-3.2) -- (gal.north);
\draw[flow] (6.1,-3.2) -- (ene.north);
\draw[flow] (10.2,-3.2) -- (col.north);

\node[tracedbox] (lin) at (2.0, -5.7)
  {\texttt{solvers.linear}\\[-1pt]{\scriptsize sparse / matrix-free}};
\node[tracedbox] (opt) at (6.1, -5.7)
  {\texttt{solvers.optimize}\\[-1pt]{\scriptsize L-BFGS, Adam}};
\node[tracedbox] (gn)  at (10.2, -5.7)
  {\texttt{Gauss--Newton}\\[-1pt]{\scriptsize least squares}};
\draw[flow] (gal) -- (lin);
\draw[flow] (ene) -- (opt);
\draw[flow] (col) -- (gn);

\node[tracedbox] (sol) at (6.1, -7.3)
  {\texttt{post}\\[-1pt]{\scriptsize Solution: evaluate, fields, norms, VTK}};
\draw[flow] (lin.south) |- (sol.west);
\draw[flow] (opt) -- (sol);
\draw[flow] (gn.south) -- (10.2,-6.5) -| (sol.north east);

\node[setupbox] (est) at (11.15, -7.3)
  {\texttt{post.estimators}\\[-1pt]{\scriptsize $\eta_e$, D\"orfler marking}};
\draw[flow] ($(sol.east)+(0,0.15)$) -- (est.west);
\draw[loop] (est.east) -- (15.15,-7.3) -- (15.15, 1.62) -- (7.4, 1.62) -- (space.north);
\node[lbl, text=blue!55!black, anchor=east, align=right] at (15.0, 1.98)
      {adapt: estimate $\to$ mark $\to$ refine\\[-2pt]
       outer Python loop --- shapes change};

\draw[back] (sol.south) -- (6.1,-8.3) -| (14.15,-2.36)
  -- (9.7,-2.36) -| ($(eval.south)+(0.8,0)$);
\node[lbl, text=red!60!black, anchor=east, align=right] at (14.0,-3.6)
      {$\nabla_{\bm{\theta}}\,,\ \nabla_{\bm{P}}$\\[-2pt]one adjoint\\[-2pt]solution};

\node[setupbox,  minimum width=5mm, minimum height=3.5mm, inner sep=1pt] at (0.1,-8.95) {};
\node[anchor=west, font=\scriptsize] at (0.42,-8.95)
      {static: shapes, connectivity, point locations};
\node[tracedbox, minimum width=5mm, minimum height=3.5mm, inner sep=1pt] at (6.6,-8.95) {};
\node[anchor=west, font=\scriptsize] at (6.92,-8.95)
      {traced: control points, weights, $\bm{\theta}$, $\bm{d}$};

\end{tikzpicture}
\caption{JAXIGA data flow. Static setup (blue) determines connectivity,
extraction and point locations. Basis evaluation (green) combines them with
traced geometry for the compute stage (orange). The dashed path attaches
first-order reverse sensitivities through one adjoint solution at a converged
regular root. The dotted adaptive loop changes array shapes in Python;
remeshing is outside the differentiation contract. Hierarchical spaces support
the two variational methods; collocation currently requires tensor-product
spaces.}
\label{fig:architecture}
\end{figure}
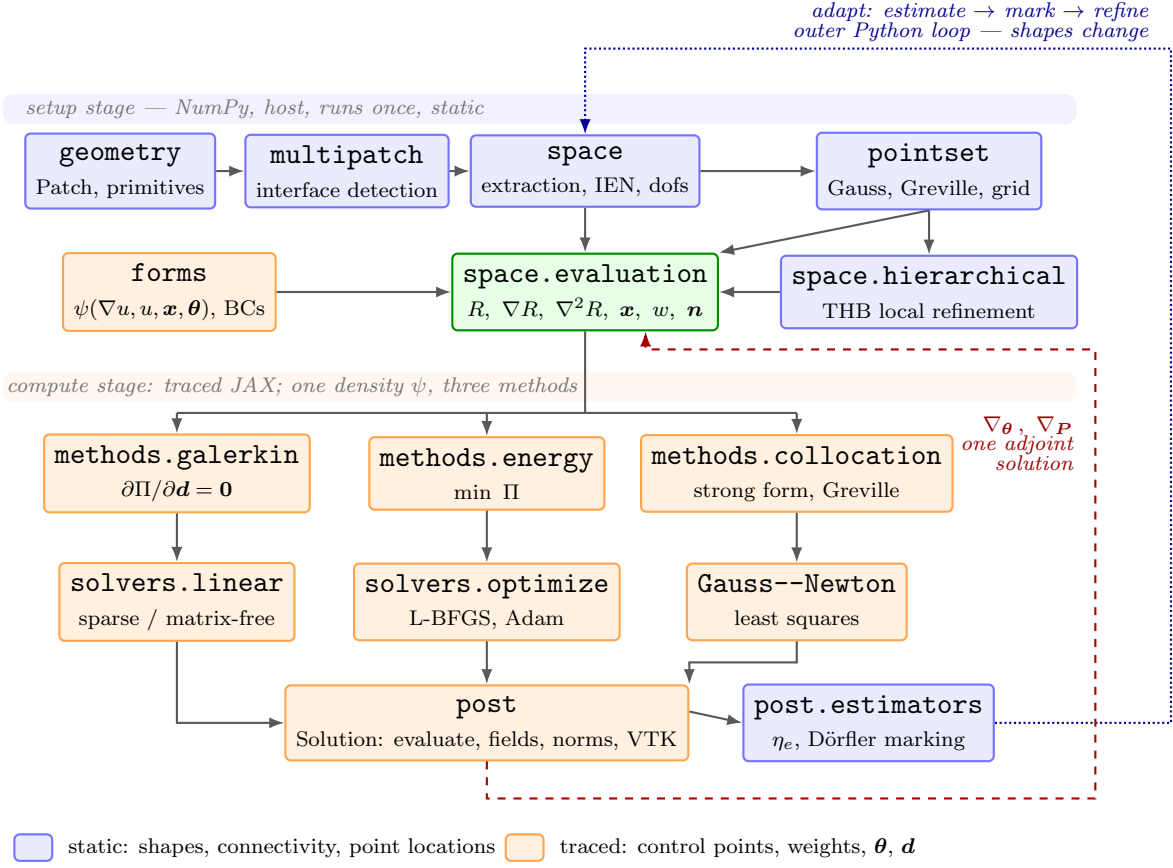

The organization of the library is shown in Figure~\ref{fig:architecture},
where the setup stage is separated from the computations traced by JAX. In
the following, we describe this separation and its consequences for assembly,
refinement and differentiation. Installation instructions and a complete
Poisson example are given in \ref{app:installation}.

\subsection{The static/traced separation}
\label{subsec:static_traced}

We divide the computation into two stages. First, the mesh is constructed on
the host CPU using NumPy, where the knot vectors, degrees, connectivity,
interface matching, boundary classification and reference point sets are
specified. The geometry, fields and material laws are then evaluated in JAX,
using traced control points, weights, extraction operators, parameters and
prescribed values. If the mesh is changed, the setup data must be rebuilt and
recompilation may be required. On the other hand, if only the parameter values
are changed, the compiled program can be reused on the same mesh.

JAX represents structured data as a \emph{pytree}, that is, a nested collection
of objects whose individual data entries are called leaves. We use this
representation to distinguish the arrays traced by JAX from the fixed data
that describe the discretization. In particular, the connectivity is stored
as immutable metadata, with a hash computed from its contents, and the knot
vectors are stored as tuples. This allows JAX to identify spaces having the
same structure and reuse the corresponding compiled program. Assembly is
vectorized over the elements and quadrature points,
while host sparse solvers are called through callbacks. Explicit derivative
rules are provided for these callbacks so that the corresponding solutions
can be differentiated.

\subsection{Geometry and multipatch coupling}
\label{subsec:design_geometry}

A patch is represented by an immutable data structure containing its knot vectors,
degrees, side labels, unweighted control points and weights. In addition to
analytically defined geometries, we can use Gmsh~\cite{gmsh2009}
\texttt{.geo}/STEP models and \texttt{.msh} meshes, or native NURBS volumes
stored as MATLAB structs. These alternatives are illustrated below:

\begin{codebox}\begin{Verbatim}[commandchars=\\\{\}]
\PY{n}{mesh} \PY{o}{=} \PY{n}{jx}\PY{o}{.}\PY{n}{read\PYZus{}gmsh}\PY{p}{(}\PY{l+s+s2}{\PYZdq{}}\PY{l+s+s2}{part.geo}\PY{l+s+s2}{\PYZdq{}}\PY{p}{,} \PY{n}{mesh\PYZus{}size}\PY{o}{=}\PY{l+s+s2}{\PYZdq{}}\PY{l+s+s2}{auto}\PY{l+s+s2}{\PYZdq{}}\PY{p}{,} \PY{n}{geometry\PYZus{}order}\PY{o}{=}\PY{l+m+mi}{3}\PY{p}{)}
\PY{n}{mesh} \PY{o}{=} \PY{n}{mesh}\PY{o}{.}\PY{n}{refine}\PY{p}{(}\PY{l+m+mi}{1}\PY{p}{)}                 \PY{c+c1}{\PYZsh{} knot insertion and CAD refit}
\PY{n}{V} \PY{o}{=} \PY{n}{jx}\PY{o}{.}\PY{n}{FunctionSpace}\PY{p}{(}\PY{n}{mesh}\PY{p}{,} \PY{n}{vec}\PY{o}{=}\PY{l+m+mi}{2}\PY{p}{)}      \PY{c+c1}{\PYZsh{} physical groups name boundaries}
\PY{c+c1}{\PYZsh{} Alternatively, import native NURBS geometry without remeshing:}
\PY{n}{patches} \PY{o}{=} \PY{n}{jx}\PY{o}{.}\PY{n}{read\PYZus{}matlab\PYZus{}nurbs}\PY{p}{(}\PY{l+s+s2}{\PYZdq{}}\PY{l+s+s2}{part.mat}\PY{l+s+s2}{\PYZdq{}}\PY{p}{)}
\end{Verbatim}
\end{codebox}

For the quadrilateral and hexahedral meshes considered here, each Gmsh cell
is converted to a tensor-product patch. The high-order polynomial map is
preserved to numerical precision by a Lagrange-to-Bernstein transformation,
and physical boundary names are retained as patch side labels for prescribing
boundary conditions. Conforming cells are then coupled with $C^0$ continuity.
For native NURBS input, the supplied knots, weights and control net are
retained, so the spline continuity within each patch is also preserved.

Knot insertion and degree elevation preserve the current geometry map.
CAD-aware enrichment instead refits boundary samples, so general CAD input
remains approximate. Exact rational reconstruction is limited to supported
conics and cylindrical extrusions. \ref{app:geometry} describes the
meshing procedure and validity checks. Differentiation starts from the
resulting fixed patch layout and does not extend through Gmsh or CAD projection.

\label{subsec:design_multipatch}
To construct a multipatch space, we first check that the control nets coincide,
the tangential knots and weights are compatible, the orientations are
admissible and the normals are antiparallel. The corresponding interface
degrees of freedom are then merged. Corner contacts in two dimensions and
edge contacts in three dimensions are classified separately, and overlaps or
nonconforming interfaces are rejected. We note that this construction gives
$C^0$ continuity across the interface, which is insufficient for primal
fourth-order models requiring $C^1$ continuity.

\subsection{Function spaces and evaluation}
\label{subsec:design_space}

A global basis is constructed from one or more patches by defining a function
space, as in the following example:

\begin{codebox}\begin{Verbatim}[commandchars=\\\{\}]
\PY{n}{V} \PY{o}{=} \PY{n}{jx}\PY{o}{.}\PY{n}{FunctionSpace}\PY{p}{(}\PY{n}{patches}\PY{p}{,} \PY{n}{vec}\PY{o}{=}\PY{l+m+mi}{2}\PY{p}{)}   \PY{c+c1}{\PYZsh{} vec: components per basis fn}
\end{Verbatim}
\end{codebox}

In this construction, B\'ezier extraction is performed on each patch, after
which the element-to-degree-of-freedom map, multipatch coupling and boundary
classification are obtained. The extraction operator for an element is the
Kronecker product of the operators in each parametric direction. We order
these operators so that the first direction varies fastest, consistently with
the ordering of the control points.

For vector-valued problems, each scalar basis function carries $n_c$ unknowns,
one for each component of the field, where $n_c$ is the value of the
\texttt{vec} argument. The global degrees of freedom are numbered one basis
function at a time, so that the components of basis function $i$ occupy the
consecutive indices $n_c i,\,n_c i+1,\ldots,n_c i+n_c-1$. In other words, with
indices starting from zero, component $c$ of basis function $i$ has the global
index $n_c i+c$.
Constrained degrees of freedom are handled by a single object that holds the
free and prescribed index sets together with the prescribed values. Its
\texttt{lift} operation reconstructs the full coefficient vector from the free
coefficients, while \texttt{restrict} extracts the free coefficients from a
full vector.
Since the prescribed values are stored in a traced array, the Dirichlet data
may depend on parameters and their contribution to the gradient is retained.

\label{subsec:design_evaluation}

The point sets used for quadrature, collocation and postprocessing have a
common representation. Their reference locations and Bernstein values are
fixed during setup, while their physical locations depend on the traced
geometry. The same evaluation procedure is used in each spatial dimension.
It applies extraction and rational correction, and then uses the geometry
Jacobian to map derivatives to physical coordinates. In this way, we obtain
the basis values, first and second derivatives, integration weights and
boundary normals. For interior Gauss rules, one reference array is shared
across the elements to avoid storing repeated data.

\subsection{Problems: one physics definition}
\label{subsec:design_problem}

A problem combines a function space and boundary conditions with a pointwise
energy density, as illustrated by the following definition:

\begin{codebox}\begin{Verbatim}[commandchars=\\\{\}]
\PY{k}{class}\PY{+w}{ }\PY{n+nc}{NeoHookean}\PY{p}{(}\PY{n}{jx}\PY{o}{.}\PY{n}{Problem}\PY{p}{)}\PY{p}{:}
    \PY{k}{def}\PY{+w}{ }\PY{n+nf}{energy}\PY{p}{(}\PY{n+nb+bp}{self}\PY{p}{,} \PY{n}{grad\PYZus{}u}\PY{p}{,} \PY{n}{u}\PY{p}{,} \PY{n}{x}\PY{p}{,} \PY{n}{p}\PY{p}{)}\PY{p}{:}
        \PY{n}{F}  \PY{o}{=} \PY{n}{jnp}\PY{o}{.}\PY{n}{eye}\PY{p}{(}\PY{l+m+mi}{2}\PY{p}{)} \PY{o}{+} \PY{n}{grad\PYZus{}u}
        \PY{n}{J}  \PY{o}{=} \PY{n}{jnp}\PY{o}{.}\PY{n}{linalg}\PY{o}{.}\PY{n}{det}\PY{p}{(}\PY{n}{F}\PY{p}{)}
        \PY{n}{I1} \PY{o}{=} \PY{n}{jnp}\PY{o}{.}\PY{n}{trace}\PY{p}{(}\PY{n}{F}\PY{o}{.}\PY{n}{T} \PY{o}{@} \PY{n}{F}\PY{p}{)}
        \PY{k}{return} \PY{p}{(}\PY{n}{p}\PY{p}{[}\PY{l+s+s2}{\PYZdq{}}\PY{l+s+s2}{mu}\PY{l+s+s2}{\PYZdq{}}\PY{p}{]} \PY{o}{/} \PY{l+m+mi}{2} \PY{o}{*} \PY{p}{(}\PY{n}{I1} \PY{o}{\PYZhy{}} \PY{l+m+mi}{2}\PY{p}{)}
                \PY{o}{\PYZhy{}} \PY{n}{p}\PY{p}{[}\PY{l+s+s2}{\PYZdq{}}\PY{l+s+s2}{mu}\PY{l+s+s2}{\PYZdq{}}\PY{p}{]} \PY{o}{*} \PY{n}{jnp}\PY{o}{.}\PY{n}{log}\PY{p}{(}\PY{n}{J}\PY{p}{)}
                \PY{o}{+} \PY{n}{p}\PY{p}{[}\PY{l+s+s2}{\PYZdq{}}\PY{l+s+s2}{lam}\PY{l+s+s2}{\PYZdq{}}\PY{p}{]} \PY{o}{/} \PY{l+m+mi}{2} \PY{o}{*} \PY{n}{jnp}\PY{o}{.}\PY{n}{log}\PY{p}{(}\PY{n}{J}\PY{p}{)} \PY{o}{*}\PY{o}{*} \PY{l+m+mi}{2}\PY{p}{)}
\end{Verbatim}
\end{codebox}

The flux and mass maps, following the JAX-FEM convention, are obtained by
differentiating this density with respect to $\nabla u$ and $u$. For
nonvariational problems, these maps may be specified directly for use with
collocation. However, an energy density is required for both variational
methods. The parameters may be supplied as a nested collection of arrays,
using the pytree representation described in
Section~\ref{subsec:static_traced}. This allows material constants and
spatially varying data to be grouped in a single argument, while derivatives
are computed with respect to the individual arrays.

To prescribe a boundary condition, one specifies a value or traction, a
boundary label or coordinate predicate, and, when needed, a component. During
setup, Dirichlet predicates select control-point indices and Neumann
predicates select boundary elements. These selections remain fixed during
differentiation, although the prescribed values remain traced. Some
coefficients, such as damage and history fields, are available only at
quadrature points. We distinguish these arrays from parameters shared by all
points using a dedicated data type. The evaluation procedure then selects
the entry for the current quadrature point and supplies it to the local
energy or constitutive function.

\subsection{Deriving three methods}
\label{subsec:one_definition}

Let $\bm{d}$ collect the free degrees of freedom and let
\begin{equation}
\Pi(\bm{d}, \bm{\theta}) \;=\; \sum_q w_q\, \psi\!\left(\nabla\bm u_q, \bm u_q, \bm{x}_q, \bm{\theta}\right)
\;-\; \sum_q w_q\, \bm{f}_q \!\cdot\! \bm{u}_q
\;-\; \sum_b \sum_q w_{bq}\, \bm{t}_{bq} \!\cdot\! \bm{u}_{bq}
\label{eq:total_energy}
\end{equation}
be the quadrature form of~\eqref{eq:energy_decomposition}, with the last sum
over Neumann boundaries. The variational drivers use~\eqref{eq:total_energy},
while collocation uses the strong form derived from $\psi$ without volume quadrature. Energy
minimization uses the optimizers of
Section~\ref{subsec:solution_strategies}.

For the Galerkin method, we solve $\partial \Pi/\partial \bm{d} = \bm{0}$.
If the energy is supplied by the user, the element tangent is obtained as the
Hessian of the element energy, with this operation vectorized over all
elements. For the built-in Poisson, linear elasticity, Kirchhoff--Love plate
and phase-field subproblems, the algebraically equivalent quadratic forms
are evaluated directly. For example, the elasticity tangent is given by
$\bm{K}_e=\sum_q w_q\bm{B}_q^{\!\top}\bm{C}_q\bm{B}_q$. If the user changes
the energy in a subclass of a built-in problem, this change is detected and
the tangent is again computed from the energy Hessian. For a linear problem,
the load vector is obtained as $-\,\mathcal{R}(\bm{0})$, where
$\mathcal{R}$ denotes the residual restricted to the free degrees of freedom.
This expression includes the source term, the Neumann work and the
contribution of inhomogeneous Dirichlet data, which therefore do not require
separate assembly. For nonlinear problems, the same tangent is used in the
Newton iteration.

The collocation driver evaluates the strong-form residual
$-\nabla\!\cdot\!\bm{\sigma} + \partial\psi/\partial u - \bm{f}$ of
Section~\ref{subsec:collocation}, with flux
$\bm\sigma=\partial\psi/\partial\nabla\bm u$. Its divergence is expanded by the chain rule
using the second physical derivatives supplied by the evaluator,
\begin{equation}
\frac{\partial \sigma_{cd}}{\partial x_d}
= \frac{\partial \sigma_{cd}}{\partial (\nabla u)_{ab}} \frac{\partial^2 u_a}{\partial x_b \partial x_d}
+ \frac{\partial \sigma_{cd}}{\partial u_a} \frac{\partial u_a}{\partial x_d}
+ \left.\frac{\partial \sigma_{cd}}{\partial x_d}\right|_{\text{explicit}} ,
\label{eq:div_flux}
\end{equation}
where repeated indices are summed and the flux derivatives are obtained by
forward-mode differentiation. A second, slower derivation is retained in which the field is
built as a function of the parametric coordinate and the divergence is taken by
nested forward-mode differentiation through it. Because this derivation does
not rely on the derivatives supplied by the evaluator, it serves as an
independent check. On multipatch
spaces, flux-continuity equations at matched interface points and Neumann
equations are appended. Finally, the resulting overdetermined system is solved by
Gauss--Newton on the normal equations.

\subsection{Solvers and implicit differentiation}
\label{subsec:design_solvers}

During setup, the library stores the unique free--free sparsity pattern in
compressed sparse row (CSR) format, together with a map from element entries
to positions in its value array. During assembly, this map is used to sum
element contributions associated with the same global matrix entry. Entries
outside the reduced free--free matrix, including padding, are excluded. The
CPU solver receives the assembled matrix values and uses the sparsity pattern
stored during setup. Matrix-free solvers apply the tangent through a residual
Jacobian--vector product.

For $\bm A\bm u=\bm b$ and a scalar loss $L$, the backward pass solves
$\bm A^{\!\top}\bm\lambda=\partial L/\partial\bm u$ and returns
\begin{equation}
 \frac{\partial L}{\partial b_i}=\lambda_i,\qquad
 \frac{\partial L}{\partial A_{ij}}=-\lambda_i u_j .
\label{eq:adjoint}
\end{equation}
Since the energy Hessian is symmetric, the same matrix values and sparsity
pattern can be used for the Galerkin adjoint problem. The accuracy of the
gradient depends on the conditioning of the problem and on the convergence
of both the primal and adjoint solutions. We note that only one adjoint system
has to be solved, independently of the number of parameters. However, the cost
of evaluating and storing the parameter derivatives may still increase with
that number.

For Newton and energy minimization, the root condition is
$\bm r=\nabla\Pi=\bm0$, while least-squares collocation uses
$\bm r=\bm J^{\!\top}\bm F=\bm0$. Once the iteration has converged to
$\bar{\bm u}$, the library detaches the iteration history from differentiation
and attaches a first-order derivative through the correction
\begin{equation}
 \bm u^*=\bar{\bm u}-\bm A(\bar{\bm u})^{-1}\bm r(\bar{\bm u}),
 \qquad \bm A=\partial\bm r/\partial\bm u .
\label{eq:implicit}
\end{equation}
At an exact regular root, this gives
$\partial\bm u^*/\partial\bm\theta=-\bm A^{-1}\bm r_{\bm\theta}$.
At a finite residual, differentiation also contributes
$\bm A^{-1}\bm A_{\bm\theta}\bm A^{-1}\bm r$. Convergence is therefore
required for the stated sensitivity formula as well as for forward accuracy.
The correction requires the solution of one additional tangent system in the
forward computation and of one adjoint system in the backward computation. Quadrature
data captured by the function remain available to automatic differentiation.
Consequently, the same expression can be used to obtain material, boundary
and geometry sensitivities.

For collocation, we need the exact derivative of the stationarity residual,
which is given by
\begin{equation}
 \frac{\partial(\bm J^{\!\top}\bm F)}{\partial\bm u}
 =\bm J^{\!\top}\bm J+\sum_i F_i\nabla^2 F_i .
\label{eq:gauss_newton_gap}
\end{equation}
We note that the second term does not necessarily vanish at the solution of an
overdetermined nonlinear least-squares problem. Therefore, although the
iteration uses Gauss--Newton, the implicit correction uses the full Hessian
of $\tfrac12\|\bm F\|^2$. In Section~\ref{subsec:ex_gradients}, we consider
a test for which omitting this term would give an incorrect sensitivity.

The differentiation results described above apply to first-order reverse-mode
derivatives computed with \texttt{grad} or \texttt{jit(grad)}, provided that
the discretization is fixed and the solution is smooth and locally
nonsingular. We refer to these requirements as the differentiation contract.

For the assembled linear solver, we provide a custom vector--Jacobian product
(VJP), which propagates derivatives of the objective backward through the
solution of the linear system. This rule supports reverse-mode
differentiation, but not the forward-mode operations \texttt{jvp} and
\texttt{jacfwd}. Moreover, repeatedly differentiating the correction
in~\eqref{eq:implicit}, with the converged iterate held independent of the
parameters, does not give a valid general rule for higher-order sensitivities.
The stated differentiation results are therefore restricted to first order.

We also check convergence within the compiled solution procedure. Newton and
optimization iterations use \texttt{lax.while\_loop}, with stopping conditions
evaluated from the current residuals. The returned solution includes residuals
and convergence flags, which must be checked before interpreting its
sensitivities. For collocation, we report both $\|\bm F\|$ and
$\|\bm J^{\!\top}\bm F\|$, since stationarity of the least-squares problem
does not require the first quantity to vanish. If an iteration fails to
converge, its result is returned without the implicit correction and is not
used for sensitivities.

\subsection{Postprocessing}
\label{subsec:design_post}

We store a computed solution together with its function space, full
coefficient vector after lifting, parameters and problem definition. These
data are used to evaluate the solution and its derivatives, compute relative
error norms, and obtain derived quantities such as stresses. The same
evaluation procedures also provide the field values needed for plots and VTK
output. In particular, a solution can be evaluated at prescribed physical
locations, as required when comparing computed fields with measurements.

To evaluate the solution at an arbitrary physical point, we first identify
the containing element by a search on the host CPU. This search provides the
element index and an initial estimate of the local coordinates. We then
invert the geometry map within the selected element using a differentiable
Newton iteration, and evaluate the spline expansion at the resulting
coordinates. The element search is a discrete operation and is not included
in automatic differentiation. Therefore, for a differentiated or
just-in-time compiled computation, the element indices and initial
coordinate estimates are prepared outside the traced computation.

During differentiation, we hold the selected element index fixed. This does
not require the element geometry or the local coordinates to remain fixed:
changes in the control points, weights or physical evaluation point can
change the inverse coordinates, and these changes propagate through the
Newton iteration and field evaluation. Changes in the solution coefficients
are also included. The resulting sensitivity applies locally, provided that
the evaluation point remains in the interior of the selected element and
the inversion has converged with a nonsingular geometry Jacobian. If a
finite update changes the containing element, the host search must be
repeated. For example, when the geometry is varied while a sensor location
is fixed in physical coordinates, the local coordinates of that sensor must
be recomputed within the selected element.

For visualization, derived fields may also be exported as cell averages.
These averages are computed by quadrature over the element interior, which
avoids evaluating singular endpoint derivatives on collapsed boundaries.

\subsection{Hierarchical refinement}
\label{subsec:design_hierarchical}

To obtain local refinement without inserting knot lines across an entire
patch, we use THB-splines~\cite{giannelli2012thb,vuong2011hierarchical}.
Related truncation-based constructions have been developed for
T-splines~\cite{wei2017truncated} and for subdivision surfaces of arbitrary
topology~\cite{wei2015thccs}.
Multilevel B\'ezier extraction~\cite{dangella2018multilevel} provides the
connectivity, extraction, control-point and weight arrays in the same form
as those used by the variational methods. The active coarse functions are
truncated in the refined regions, while the partition of unity and the
original geometry are preserved. Since the number of locally supported
functions can vary between elements, zero rows are added to the extraction
arrays to allow vectorized evaluation.

Refinement is mirrored across patch interfaces so that the resulting space
remains conforming. We also impose class-2
admissibility~\cite{buffa2016adaptive}, which limits the overlap of refinement
levels over the supports of the basis functions. This gives stronger control
over the number of active functions on an element than a restriction on
neighboring element sizes alone.

\begin{codebox}\begin{Verbatim}[commandchars=\\\{\}]
\PY{k+kn}{import}\PY{+w}{ }\PY{n+nn}{jaxiga}\PY{+w}{ }\PY{k}{as}\PY{+w}{ }\PY{n+nn}{jx}

\PY{k}{def}\PY{+w}{ }\PY{n+nf}{make\PYZus{}problem}\PY{p}{(}\PY{n}{V}\PY{p}{)}\PY{p}{:}
\PY{+w}{    }\PY{l+s+sd}{\PYZdq{}\PYZdq{}\PYZdq{}A factory: BCs resolve against whichever space is current.\PYZdq{}\PYZdq{}\PYZdq{}}
    \PY{k}{return} \PY{n}{jx}\PY{o}{.}\PY{n}{Poisson}\PY{p}{(}\PY{n}{V}\PY{p}{,} \PY{n}{dirichlet}\PY{o}{=}\PY{p}{[}\PY{n}{jx}\PY{o}{.}\PY{n}{DirichletBC}\PY{p}{(}\PY{n}{g}\PY{p}{,} \PY{n}{where}\PY{o}{=}\PY{n}{on\PYZus{}bnd}\PY{p}{)}\PY{p}{]}\PY{p}{)}

\PY{n}{V} \PY{o}{=} \PY{n}{jx}\PY{o}{.}\PY{n}{FunctionSpace}\PY{p}{(}\PY{n}{patches}\PY{p}{)}
\PY{k}{for} \PY{n}{cycle} \PY{o+ow}{in} \PY{n+nb}{range}\PY{p}{(}\PY{n}{n\PYZus{}cycles}\PY{p}{)}\PY{p}{:}
    \PY{n}{sol} \PY{o}{=} \PY{n}{jx}\PY{o}{.}\PY{n}{solve}\PY{p}{(}\PY{n}{make\PYZus{}problem}\PY{p}{(}\PY{n}{V}\PY{p}{)}\PY{p}{,} \PY{n}{params}\PY{o}{=}\PY{n}{params}\PY{p}{)}
    \PY{n}{eta} \PY{o}{=} \PY{n}{jx}\PY{o}{.}\PY{n}{residual\PYZus{}indicator}\PY{p}{(}\PY{n}{sol}\PY{p}{,} \PY{n}{params}\PY{p}{)}
    \PY{n}{V} \PY{o}{=} \PY{n}{jx}\PY{o}{.}\PY{n}{refine\PYZus{}elements}\PY{p}{(}\PY{n}{V}\PY{p}{,} \PY{n}{jx}\PY{o}{.}\PY{n}{dorfler\PYZus{}mark}\PY{p}{(}\PY{n}{eta}\PY{p}{,} \PY{n}{frac}\PY{o}{=}\PY{l+m+mf}{0.3}\PY{p}{)}\PY{p}{)}

\PY{c+c1}{\PYZsh{} the same loop, as a driver}
\PY{n}{sol}\PY{p}{,} \PY{n}{V}\PY{p}{,} \PY{n}{history} \PY{o}{=} \PY{n}{jx}\PY{o}{.}\PY{n}{adapt}\PY{p}{(}\PY{n}{make\PYZus{}problem}\PY{p}{,} \PY{n}{params}\PY{p}{,} \PY{n}{space}\PY{o}{=}\PY{n}{V0}\PY{p}{,}
                           \PY{n}{n\PYZus{}cycles}\PY{o}{=}\PY{l+m+mi}{8}\PY{p}{,} \PY{n}{frac}\PY{o}{=}\PY{l+m+mf}{0.3}\PY{p}{)}
\end{Verbatim}
\end{codebox}

Marking, that is, the selection of the elements to be refined on the basis of
an error indicator, and the subsequent refinement both change array shapes, so
adaptivity is implemented in an outer Python loop. Gradients do not propagate across refinement cycles, but
computations on each fixed adapted space remain covered by the differentiation
contract. The marking indicator
$\eta_e^2=h_e^2\sum_q w_q|\mathcal R(u_h)|^2$ reuses the collocation strong
residual at Gauss points. Interior flux jumps vanish for smooth coefficients
and a $C^1$ basis, but repeated knots, patch interfaces, coefficient
discontinuities and boundary residuals require additional terms that this
volume-only indicator omits. Its reliability therefore does not extend to
general problems with these features.

When the mesh is refined, the state variables must also be transferred to the
new space. For nested spline fields, we use $L^2$ projection, which reproduces
the coarse field provided that the quadrature is consistent and the discrete
mass matrix is nonsingular. The fracture history is a running maximum stored
at quadrature points. For this quantity, we use an approximate transfer from
the nearest parent point in order to avoid overshoots due to projection.

\subsection{Fourth-order problems}
\label{subsec:design_hessian}

We next consider primal fourth-order formulations, for which the required
continuity is available within a patch when maximally smooth splines of
degree $p\geq2$ are used. These splines provide $C^1$ continuity. By setting
\texttt{needs\_hessian}, a problem can request the physical second derivatives
in its energy kernel. For general second-gradient energies, the element
Hessians are computed by automatic differentiation. For the built-in plate
problem, they are obtained by contracting the curvature--displacement matrix
with the bending tensor.

The Kirchhoff--Love shell uses the same flag on a surface patch whose control
points lie in $\mathbb{R}^3$. For such a patch, we retain derivatives with
respect to the two parametric coordinates rather than transform them to
Cartesian derivatives using the volume-map procedure. The evaluator also
computes the surface measure $|\bm{a}_1\times\bm{a}_2|$. The shell energy
receives $\bm{u}_{,\alpha}$ and $\bm{u}_{,\alpha\beta}$ together with the
reference metric at each quadrature point, so that the curved geometry is
accounted for in the strain expressions. We express the strains in terms of
the deformed surface and use
forward-mode differentiation to obtain their exact linearization, as shown
in the following example:

\begin{codebox}\begin{Verbatim}[commandchars=\\\{\}]
\PY{k}{class}\PY{+w}{ }\PY{n+nc}{KirchhoffLoveShell}\PY{p}{(}\PY{n}{jx}\PY{o}{.}\PY{n}{Problem}\PY{p}{)}\PY{p}{:}
    \PY{n}{needs\PYZus{}hessian} \PY{o}{=} \PY{k+kc}{True}    \PY{c+c1}{\PYZsh{} kernel receives parametric u\PYZus{},a and u\PYZus{},ab}

    \PY{k}{def}\PY{+w}{ }\PY{n+nf}{strains}\PY{p}{(}\PY{n+nb+bp}{self}\PY{p}{,} \PY{n}{grad\PYZus{}u}\PY{p}{,} \PY{n}{hess\PYZus{}u}\PY{p}{,} \PY{n}{p}\PY{p}{)}\PY{p}{:}
        \PY{k}{def}\PY{+w}{ }\PY{n+nf}{local}\PY{p}{(}\PY{n}{du}\PY{p}{,} \PY{n}{d2u}\PY{p}{)}\PY{p}{:}    \PY{c+c1}{\PYZsh{} exact strains of the deformed midsurface}
            \PY{n}{E}\PY{p}{,} \PY{n}{K} \PY{o}{=} \PY{n+nb+bp}{self}\PY{o}{.}\PY{n}{\PYZus{}covariant\PYZus{}strains}\PY{p}{(}\PY{n}{du}\PY{p}{,} \PY{n}{d2u}\PY{p}{,} \PY{n}{p}\PY{p}{[}\PY{l+s+s2}{\PYZdq{}}\PY{l+s+s2}{a}\PY{l+s+s2}{\PYZdq{}}\PY{p}{]}\PY{p}{,} \PY{n}{p}\PY{p}{[}\PY{l+s+s2}{\PYZdq{}}\PY{l+s+s2}{da}\PY{l+s+s2}{\PYZdq{}}\PY{p}{]}\PY{p}{)}
            \PY{k}{return} \PY{n+nb+bp}{self}\PY{o}{.}\PY{n}{\PYZus{}to\PYZus{}local}\PY{p}{(}\PY{n}{E}\PY{p}{,} \PY{n}{p}\PY{p}{[}\PY{l+s+s2}{\PYZdq{}}\PY{l+s+s2}{a}\PY{l+s+s2}{\PYZdq{}}\PY{p}{]}\PY{p}{)}\PY{p}{,} \PY{n+nb+bp}{self}\PY{o}{.}\PY{n}{\PYZus{}to\PYZus{}local}\PY{p}{(}\PY{n}{K}\PY{p}{,} \PY{n}{p}\PY{p}{[}\PY{l+s+s2}{\PYZdq{}}\PY{l+s+s2}{a}\PY{l+s+s2}{\PYZdq{}}\PY{p}{]}\PY{p}{)}
        \PY{n}{zero} \PY{o}{=} \PY{p}{(}\PY{n}{jnp}\PY{o}{.}\PY{n}{zeros\PYZus{}like}\PY{p}{(}\PY{n}{grad\PYZus{}u}\PY{p}{)}\PY{p}{,} \PY{n}{jnp}\PY{o}{.}\PY{n}{zeros\PYZus{}like}\PY{p}{(}\PY{n}{hess\PYZus{}u}\PY{p}{)}\PY{p}{)}
        \PY{c+c1}{\PYZsh{} jvp returns (value, tangent); the value at u = 0 is not needed}
        \PY{n}{\PYZus{}}\PY{p}{,} \PY{p}{(}\PY{n}{eps}\PY{p}{,} \PY{n}{kappa}\PY{p}{)} \PY{o}{=} \PY{n}{jax}\PY{o}{.}\PY{n}{jvp}\PY{p}{(}\PY{n}{local}\PY{p}{,} \PY{n}{zero}\PY{p}{,} \PY{p}{(}\PY{n}{grad\PYZus{}u}\PY{p}{,} \PY{n}{hess\PYZus{}u}\PY{p}{)}\PY{p}{)}
        \PY{k}{return} \PY{n}{eps}\PY{p}{,} \PY{n}{kappa}      \PY{c+c1}{\PYZsh{} exact linearization by forward\PYZhy{}mode AD}

    \PY{k}{def}\PY{+w}{ }\PY{n+nf}{energy}\PY{p}{(}\PY{n+nb+bp}{self}\PY{p}{,} \PY{n}{grad\PYZus{}u}\PY{p}{,} \PY{n}{u}\PY{p}{,} \PY{n}{x}\PY{p}{,} \PY{n}{p}\PY{p}{,} \PY{n}{hess\PYZus{}u}\PY{o}{=}\PY{k+kc}{None}\PY{p}{)}\PY{p}{:}
        \PY{n}{eps}\PY{p}{,} \PY{n}{kappa} \PY{o}{=} \PY{n+nb+bp}{self}\PY{o}{.}\PY{n}{strains}\PY{p}{(}\PY{n}{grad\PYZus{}u}\PY{p}{,} \PY{n}{hess\PYZus{}u}\PY{p}{,} \PY{n}{p}\PY{p}{)}
        \PY{n}{D}\PY{p}{,} \PY{n}{t} \PY{o}{=} \PY{n+nb+bp}{self}\PY{o}{.}\PY{n}{material}\PY{p}{(}\PY{n}{p}\PY{p}{)}\PY{p}{,} \PY{n}{p}\PY{p}{[}\PY{l+s+s2}{\PYZdq{}}\PY{l+s+s2}{t}\PY{l+s+s2}{\PYZdq{}}\PY{p}{]}
        \PY{k}{return} \PY{l+m+mf}{0.5} \PY{o}{*} \PY{n}{t} \PY{o}{*} \PY{n}{eps} \PY{o}{@} \PY{n}{D} \PY{o}{@} \PY{n}{eps} \PY{o}{+} \PY{n}{t}\PY{o}{*}\PY{o}{*}\PY{l+m+mi}{3} \PY{o}{/} \PY{l+m+mi}{24} \PY{o}{*} \PY{n}{kappa} \PY{o}{@} \PY{n}{D} \PY{o}{@} \PY{n}{kappa}

\PY{n}{V} \PY{o}{=} \PY{n}{jx}\PY{o}{.}\PY{n}{FunctionSpace}\PY{p}{(}\PY{n}{roof\PYZus{}surface}\PY{p}{,} \PY{n}{vec}\PY{o}{=}\PY{l+m+mi}{3}\PY{p}{)}   \PY{c+c1}{\PYZsh{} 2D patch, 3D control points}
\PY{n}{problem} \PY{o}{=} \PY{n}{KirchhoffLoveShell}\PY{p}{(}\PY{n}{V}\PY{p}{,} \PY{n}{source}\PY{o}{=}\PY{n}{self\PYZus{}weight}\PY{p}{,} \PY{n}{dirichlet}\PY{o}{=}\PY{n}{diaphragms}\PY{p}{)}
\PY{n}{solution} \PY{o}{=} \PY{n}{jx}\PY{o}{.}\PY{n}{solve}\PY{p}{(}\PY{n}{problem}\PY{p}{,} \PY{n}{params}\PY{o}{=}\PY{n}{problem}\PY{o}{.}\PY{n}{parameters}\PY{p}{(}\PY{n}{E}\PY{p}{,} \PY{n}{nu}\PY{p}{,} \PY{n}{t}\PY{p}{)}\PY{p}{)}
\end{Verbatim}
\end{codebox}

A clamped edge requires additional constraints because fixing the normal slope
of a spline means constraining two layers of control points rather than one.
A simply supported edge or a shell diaphragm can use the existing boundary
condition procedures, since the zero-moment condition follows naturally from
the bending energy.

\subsection{Assembly variants, batching and time integration}
\label{subsec:design_variants}

The tangent can also be applied to a vector through the residual Jacobian,
without assembling a global matrix. In this matrix-free formulation, the
adjoint is computed using \texttt{jax.lax.custom\_linear\_solve} with the
transposed operator. The Krylov vectors require $O(n_{\mathrm{dofs}})$
storage, although additional storage is needed for the basis arrays and the
local tangents used in Jacobi preconditioning. To reduce the storage of basis
data, assembly can be performed in chunks. This replaces the
$O(n_e n_q n_{\mathrm{loc}}(1+d))$ storage for all elements by the storage
needed for one chunk. The basis data can also be recomputed during the
backward pass through checkpointing, reducing the amount of data retained.

For a fixed sparsity pattern, several sets of matrix values and right-hand
sides can be processed using \texttt{jax.vmap}. This batching procedure has
been tested for all linear solver backends, although host (CPU) callbacks process
the members of a batch sequentially. We have not included nonlinear batching
in these tests. For time-dependent problems, the Newmark-$\beta$ and
generalized-$\alpha$ methods~\cite{chung1993generalized} are implemented using
\texttt{jax.lax.scan}, whose reverse pass differentiates the trajectory. If
the mass, stiffness and time-step parameters are constant, the effective
matrix is formed once before the time loop. Numerical checks of these
procedures are given in Section~\ref{subsec:ex_variants}, together with an
application to full-waveform inversion~\cite{virieux2009fwi}.

\section{Numerical Examples}
\label{sec:examples}

In this section, we present numerical experiments to examine the accuracy of
the methods, the computation of sensitivities and the cost of the
implementation. We compare the computed fields with analytical solutions
and independent implementations, and verify the sensitivities using
finite-difference approximations.

For the CPU and GPU plate timings, we use the same benchmark machine,
specified in Table~\ref{tab:provenance}. The other examples were measured
separately, so their timings are used for comparisons within the respective
studies rather than across machines or software versions. CPU measurements
described as ``single-core'' follow a single-core measurement protocol;
processor affinity was not verified for every computation. The cost of
initial compilation and that of repeated solutions are reported separately
where appropriate.

\begin{table}[htbp]
\centering
\caption{Computational environments. The benchmark CPU and GPU rows identify the machine used for the plate timing
comparison, and the final row gives the material-identification hardware.
The software versions describe the main result-generation environment.
The hardware for the fracture timings is not fully specified; these timings
are used only for comparisons within the fracture examples.}

\label{tab:provenance}
\small
\begin{tabular}{l l}
\hline
JAXIGA & 0.3.0  \\
JAX / jaxlib & 0.11.1 / 0.11.1 \\
NumPy / SciPy & 2.5.2 / 1.18.0 \\
DOLFINx / basix / UFL & 0.9.0 / 0.9.0 / 2024.2.0 (Python 3.12.13) \\
Python & 3.14.7 \\
Platform & macOS-26.6.2-arm64-arm-64bit-Mach-O \\
Backend, device & cpu \\
Precision & float64 \\
Random seed & 20260901 \\
Benchmark machine, CPU & AMD EPYC 7J13 64-Core Processor \\
\quad accelerator & NVIDIA A100-SXM4-40GB (CUDA) \\
Material identification (Section~\ref{subsec:ex_mnist}) & Apple M1 Pro CPU; NVIDIA A10 (CUDA 12) \\
\hline
\end{tabular}
\end{table}

\subsection{Convergence of the Galerkin discretization}
\label{subsec:ex_convergence}

\begin{figure}[htbp]
\centering
\includegraphics[width=\textwidth]{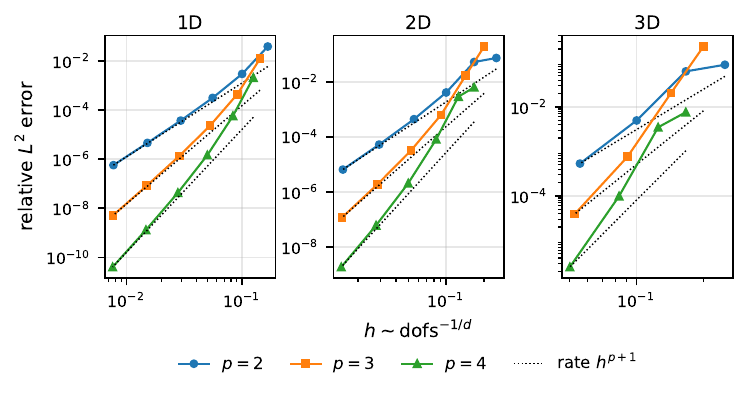}
\caption{Relative $L^2$ error against mesh size for the Poisson problem in one,
two and three dimensions, over six refinement levels in one and two dimensions
and four in three. Dotted lines indicate the theoretical rate $h^{p+1}$,
anchored at the finest mesh. The coarsest meshes resolve fewer than two
elements per wavelength, which is where the higher degrees are visibly
pre-asymptotic and the $p=4$ curve crosses the others.}
\label{fig:poisson_convergence}
\end{figure}

We first consider $-\Delta u = f$ with homogeneous Dirichlet data on the unit
interval, square and cube. The right-hand side $f$ is chosen so that the exact
solution is $u = \prod_{i=1}^{d}\sin(2\pi x_i)$. We use six uniform refinement
levels in one and two dimensions, and four levels in three dimensions, where
the storage required for the assembled tangent limits further refinement.
The errors are shown in Figure~\ref{fig:poisson_convergence}, together with
the theoretical convergence rate. On the coarsest meshes, the solution is
under-resolved and the errors have not yet reached the asymptotic regime.

\begin{table}[htbp]
\centering
\caption{Observed asymptotic convergence rates for the Poisson problem, taken
from the finest two refinement levels in each case.}
\label{tab:poisson_convergence}
\small
\begin{tabular}{c c c c c}
\hline
& \multicolumn{3}{c}{observed rate} & \\
$p$ & $d=1$ & $d=2$ & $d=3$ & theory \\
\hline
2 & 3.00 & 3.02 & 3.22 & 3 \\
3 & 4.00 & 4.02 & 4.27 & 4 \\
4 & 5.00 & 5.02 & 5.29 & 5 \\
\hline
\end{tabular}
\end{table}

From Table~\ref{tab:poisson_convergence}, we observe that the convergence
rates approach $p+1$ in one and two dimensions. On the finest one-dimensional
meshes, the rates are $3.00$, $4.00$ and $5.00$ for $p = 2, 3, 4$,
respectively, while the corresponding two-dimensional rates are $3.02$,
$4.02$ and $5.02$. In three dimensions, the rates remain slightly above the
asymptotic values, since four refinement levels are not sufficient to reach
the asymptotic regime. They agree with the two-dimensional rates at the same
resolution per wavelength. We use the same problem class, solver and error
norms in all dimensions, changing only the geometry constructor.

\subsection{Imported CAD and NURBS geometries}
\label{subsec:ex_geometry}

\begin{figure}[htbp]
\centering
\includegraphics[width=\textwidth]{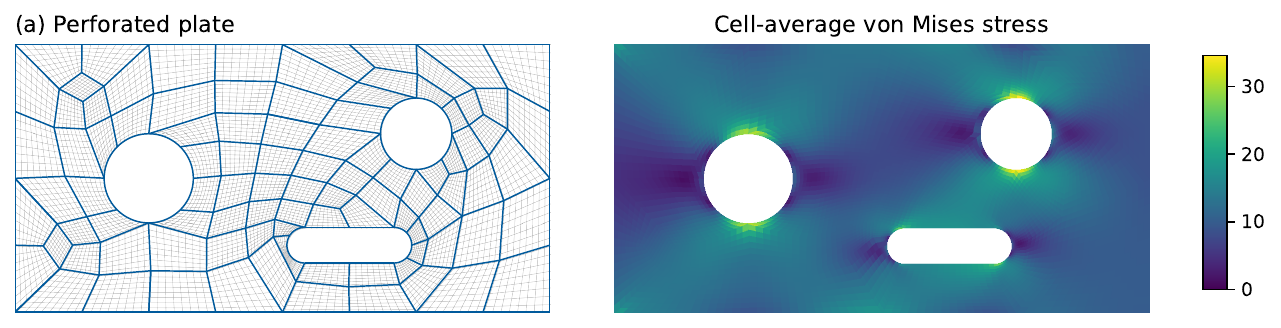}
\includegraphics[width=\textwidth]{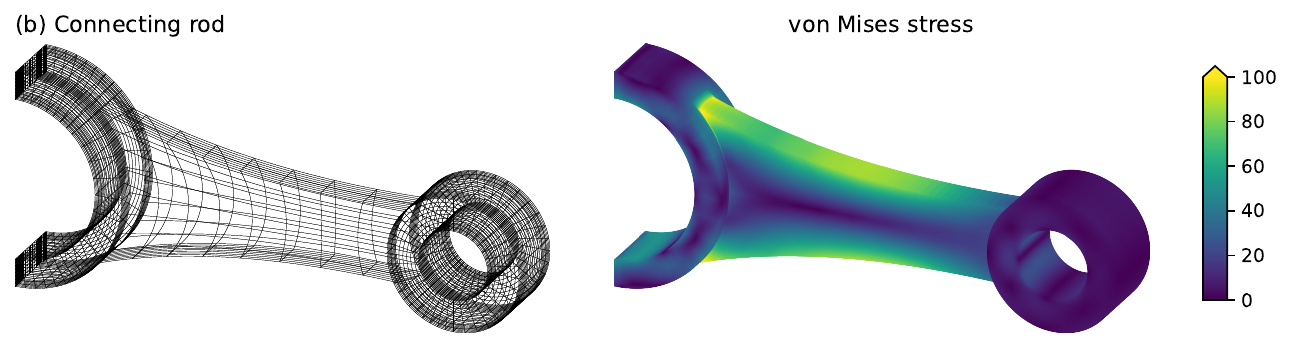}
\caption{Imported geometries with exclusively quad/hex analysis meshes.
The left panels show the undeformed meshes, with thick blue lines in (a)
marking the 95 initial patches. The right panels show the von Mises stress on
separate color scales. For the plate, area-weighted cell averages are plotted.
For the rod, the stress is evaluated pointwise on the boundary faces, so that it
varies inside each element. Its color scale is capped at 100, a value exceeded
by $0.1\%$ of the surface samples, which lie near the re-entrant edges where the
shank joins the big end (maximum 156). The plate has
6,080 quadrilateral elements and 19,496 displacement degrees of freedom, while
the rod has 12,800 hexahedral knot spans and 63,720 degrees of freedom. The plate
uses a CAD-refitted polynomial map, and the rod retains its native rational geometry.}
\label{fig:imported_geometry}
\end{figure}

In Figure~\ref{fig:imported_geometry}, we consider two examples of imported
geometry for linear elasticity. The first is a $12\times6$ plate with unequal circular holes and an offset
rounded slot, which is meshed from a Gmsh \texttt{.geo} model using automatic coarse
sizing, cubic geometry and three uniform refinement levels with CAD refitting.
Physical groups
identify the fixed left edge and right-edge traction $(10,0)$, while the other edges
are traction-free. We use plane stress with $E=10^5$ and $\nu=0.3$.
The sampled CAD boundary distance decreases from $3.16\times10^{-3}$ to
$1.35\times10^{-5}$ in model units after refitting.

For the second example, we import the 25 native NURBS volumes of a connecting
rod supplied with
\href{https://github.com/canitesc/IGAPack}{IGAPack}'s
\texttt{ConnectingRodMP.m}. After orientation correction, quadratic degree
elevation and three uniform knot-refinement levels, 39 patch interfaces are coupled
conformingly. The big-end cut faces are fixed and the small-end bore carries
traction $(0,0,-1)$, with $E=4\times10^5$, $\nu=0.3$, in consistent model units.
Both problems are solved using Jacobi-preconditioned CG with a relative
tolerance of $10^{-8}$.
The standalone scripts also export displacement and stress to VTK.

\subsection{Agreement between the three methods}
\label{subsec:ex_methods}

We next compare the three numerical methods on the two-dimensional Poisson
problem. Since the Galerkin residual is the gradient of the discrete energy,
the Galerkin and energy-minimization methods have the same solution for a
strictly convex problem, provided that both iterations converge. The results
are given in Table~\ref{tab:method_comparison}, where we also report the
collocation errors on the same meshes.

\begin{table}[htbp]
\centering
\caption{The three solution methods on the two-dimensional Poisson problem.
The method columns report relative $L^2$ errors.
Galerkin and energy minimization agree to solver tolerance, as shown by the
relative difference in the last column. Collocation is a different discretization
and converges at a lower rate, as theory predicts.}
\label{tab:method_comparison}
\small
\begin{tabular}{c c r l l l l}
\hline
$p$ & ref. & dofs & Galerkin & energy min. & collocation & $\|u_{\mathrm{EM}}-u_{\mathrm{G}}\|_\infty/\|u\|_\infty$ \\
\hline
2 & 3 & 100 & 4.067e-03 & 4.067e-03 & 2.721e-02 & 1.0e-15 \\
2 & 4 & 324 & 4.362e-04 & 4.362e-04 & 6.515e-03 & 8.9e-16 \\
2 & 5 & 1,156 & 5.226e-05 & 5.226e-05 & 1.612e-03 & 1.4e-15 \\
3 & 3 & 121 & 6.182e-04 & 6.182e-04 & 4.879e-02 & 7.2e-16 \\
3 & 4 & 361 & 3.204e-05 & 3.204e-05 & 1.275e-02 & 1.2e-15 \\
3 & 5 & 1,225 & 1.899e-06 & 1.899e-06 & 3.208e-03 & 2.2e-15 \\
4 & 3 & 144 & 8.067e-05 & 8.067e-05 & 1.002e-03 & 4.8e-15 \\
4 & 4 & 400 & 2.060e-06 & 2.060e-06 & 7.370e-05 & 3.0e-15 \\
4 & 5 & 1,296 & 6.051e-08 & 6.051e-08 & 4.983e-06 & 5.1e-15 \\
\hline
\end{tabular}
\end{table}

For all nine cases, the relative difference between the Galerkin and
energy-minimization solutions is at most $5.1 \times 10^{-15}$, which is of
the order of round-off error. The two implementations share only the energy
definition and its evaluation, so this agreement provides a check of the
assembly, the imposition of boundary conditions and the linear solver. We
also verify that the Galerkin residual vanishes at the energy minimizer. For
a hyperelastic problem, the solutions obtained by Newton iteration and energy
minimization agree to a relative difference of $1.5\times10^{-13}$, and the
Newton iteration
has the expected quadratic convergence.

The collocation solution is less accurate on the same mesh and agrees with
the variational solutions only up to discretization error, since it is
obtained from a different discrete problem. We observe rates of approximately
$p-1$ for odd degree and $p$ for even degree, consistently with the classical
convergence result for collocation at Greville abscissae. This difference
between odd and even degrees provides a check of the strong-form residual.
We also compare the two independent ways of computing this residual, one
using second physical derivatives and the other using nested differentiation
of the parametric field. On a curved NURBS annulus, the two results agree to
a relative difference of $10^{-15}$. This example requires the full chain rule, including
the transformation of second derivatives to physical coordinates.

\subsection{Verification against an independent implementation}
\label{subsec:ex_parity}

We compare JAXIGA with independent FEniCSx~\cite{dolfinx,basix} and
JAX-FEM~\cite{jaxfem} implementations of the Kirsch plate problem. The domain
is a square of half-width four with a unit circular hole. The hole is
traction-free, and outer tractions come from the analytical infinite-plate
solution. Prescribing $u_y$ at $(\pm1,0)$ and $u_x$ at $(0,\pm1)$ removes
rigid modes consistently with that solution. The $L^2$ error calculation also removes the
best-fitting rigid motion.

\begin{figure}[!htbp]
\centering
\includegraphics[width=\textwidth]{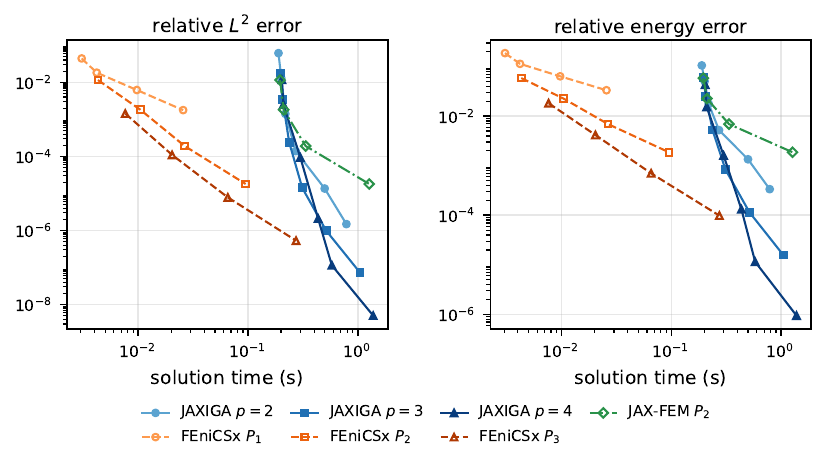}
\caption{Plate with a circular hole: relative $L^2$ error (left) and energy
error (right) against the median time of repeated solutions in double precision. JAXIGA
($p=2,3,4$, solid) and JAX-FEM ($P_2$, dash--dot) run on the GPU. FEniCSx
($P_1,P_2,P_3$, dashed) is the CPU reference from the same machine
(Table~\ref{tab:provenance}). Timings include assembly, boundary conditions and
the solution of the linear system, but exclude mesh construction and initial compilation.
JAXIGA uses the solver selection in Section~\ref{subsec:solution_strategies}, while
JAX-FEM uses SuperLU and FEniCSx uses MUMPS. The JAXIGA timings use individual
calls to its operations, without compiling the entire solution procedure as
one function. The resulting overhead for small systems is discussed in the text.}
\label{fig:plate_work_precision}
\end{figure}

JAXIGA uses four exact NURBS patches cut along the diagonals, each joining a
$90^\circ$ hole arc to a square side. FEniCSx uses $P_1$--$P_3$ fields on
curved quadratic triangles generated by an independent transfinite map.
JAX-FEM uses the same geometry and connectivity as FEniCSx $P_2$, with
\texttt{TRI6} elements and eighth-order quadrature. The reference solvers are
PETSc/MUMPS for FEniCSx and SciPy SuperLU for JAX-FEM. The reference
implementations use separate environments with fixed software versions,
and their numerical results are exchanged as JSON files.

As shown in Figure~\ref{fig:plate_work_precision}, all three codes have the
expected convergence rates. The JAX-FEM and FEniCSx $P_2$
results agree to within relative differences of $9.9\times10^{-9}$ in the
$L^2$ error, $2.5\times10^{-13}$ in the energy error and $4.9\times10^{-14}$
in total strain energy. We also observe that JAXIGA requires fewer unknowns for a
comparable error. In particular, degree four gives a relative error of
$7.8\times10^{-5}$ in the $L^2$ norm with $2{,}112$ unknowns, while FEniCSx
$P_3$ gives $1.1\times10^{-4}$ with $4{,}800$ unknowns.

However, a smaller number of unknowns does not necessarily imply a lower
computational cost at every level of accuracy. On the CPU,
FEniCSx $P_3$ reaches $7.8\times10^{-6}$ error in $0.065$\,s and
$5.3\times10^{-7}$ in $0.27$\,s, while JAXIGA degree four reaches
$2.1\times10^{-6}$ in $0.38$\,s. For JAXIGA on the NVIDIA A100 GPU,
the measurements use individual calls to the assembly and solution
operations, rather than compiling the entire procedure as one function.
With this execution procedure, the overhead for small systems is
approximately $0.19$\,s, and
the finest discretization gives an error of
$5.0\times10^{-9}$ in $1.38$\,s. At an intermediate resolution, JAXIGA
requires $0.43$\,s against JAX-FEM's $1.27$\,s and is nearly nine times more
accurate. Here we compare the specified public solver interfaces, including
the residual and Jacobian checks performed by JAX-FEM after each update.

In Table~\ref{tab:plate_cost}, we examine the effect of the solver at the
finest JAXIGA resolution. Although assembly is accelerated by $24\times$ on
the A100, the total cost is still dominated by the host SuperLU solver. When
CG is performed on the device, the total time is reduced to $1.38$\,s, and
the result agrees with the direct solution to a relative difference of
$5\times10^{-11}$.
The assembly map constructed during setup combines $5.12$ million local entries into
$813{,}364$ unique free--free nonzeros. The solver selection described in
Section~\ref{subsec:solution_strategies} chooses the faster solver on 17 of
the 18 plotted meshes. For the remaining mesh, the selected solver takes
$11\%$ longer. We use host SuperLU for the JAX-FEM curves, since its device
BiCGSTAB solver fails the convergence assertion from $2{,}048$ cells on this
benchmark.

\begin{table}[htbp]
\centering
\caption{JAXIGA plate benchmark at $36{,}448$ unknowns. Assembly and total
times for repeated solutions on the AMD EPYC CPU and NVIDIA A100 GPU, excluding initial compilation.
For this benchmark, the reduction in assembly time gives a substantial
reduction in total time only when the linear system is also solved on the device.}
\label{tab:plate_cost}
\small
\begin{tabular}{l l r r}
\toprule
Backend & Linear solver & Assembly (s) & Total (s) \\
\midrule
CPU & SuperLU & 1.46 & 31.8 \\
GPU & Host SuperLU & 0.060 & 29.9 \\
GPU & Device CG & 0.060 & 1.38 \\
\bottomrule
\end{tabular}
\end{table}

We also compare the computed fields with a FEniCSx reference having
$167{,}040$ degrees of freedom. This reference differs from the analytical
solution by $2.0\times10^{-7}$ at 410 interior points. The JAXIGA errors
measured against the reference and against the analytical solution agree to
three to seven significant digits. On the finest meshes, the
total strain energies are $0.0337796503$ (JAXIGA) and $0.0337796508$
(FEniCSx), against the analytical value $0.0337796508$. Integrating the
analytical energy with the two geometries and quadratures agrees to a relative
difference of $5.3\times10^{-10}$.

\begin{figure}[htbp]
\centering
\includegraphics[width=\textwidth]{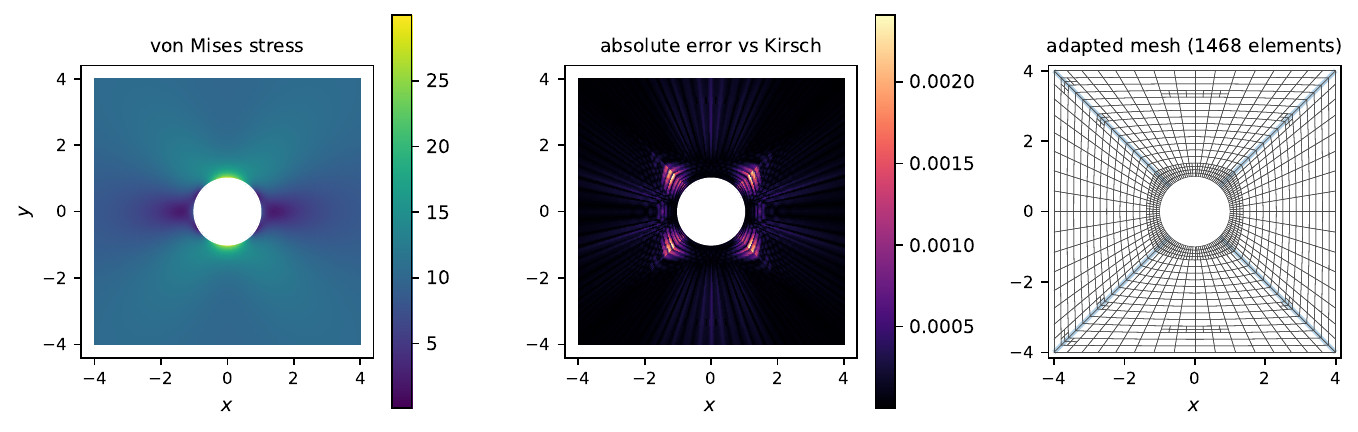}
\caption{Plate with a circular hole, computed with degree four on the
adaptively refined mesh. The hole is excluded from the plotted domain. The
left panel shows the von Mises stress, with the expected concentration factor
of three. The center panel gives the absolute error relative to the
analytical Kirsch solution, which is concentrated near the hole and along
the four patch interfaces. The right panel shows the adapted mesh, with
interfaces marked in blue. Refinement is mirrored across these interfaces,
so the discretization remains conforming across the diagonals.}
\label{fig:plate}
\end{figure}

Figure~\ref{fig:plate} shows the von Mises stress, its error against the Kirsch
solution and an adapted mesh. The expected stress concentration factor of three
is recovered, and the error concentrates near the hole and along the four $C^0$
interfaces. Starting from the coarsest mesh, residual-driven refinement with
D\"orfler fraction $0.5$ reaches an energy-norm error of $3.4\times10^{-5}$ with
$3{,}856$ degrees of freedom on $1{,}468$ elements over five levels, whereas the
uniform family needs about $5{,}400$. This factor of $1.4$ is a constant gain only,
since the analytic Kirsch solution already converges at the optimal uniform
rate $O(N^{-p/2})$. For the corner singularity in
Section~\ref{subsec:ex_adaptive}, adaptivity also improves the convergence
rate. The volume-residual indicator does not measure flux
jumps across the $C^0$ interfaces, which remains a future improvement.

To complement these comparisons of complete solutions, we verify individual
computational procedures against a pure-NumPy implementation of the same
discretization, using a relative tolerance of $10^{-12}$. These checks cover
B\'ezier extraction for degrees two to five, including repeated interior knots,
connectivity and boundary classification, control nets after degree elevation
and knot insertion, and basis values and derivatives at Gauss points. The
evaluated basis data agree with the NumPy reference to between
$2\times10^{-16}$ and $4\times10^{-15}$. The reference implementation uses
explicit loops and is intended for these individual checks, rather than for
the larger numerical examples.

We also verify mathematical identities and geometric properties without
using a reference implementation. These checks include partition of unity,
vanishing sums of basis derivatives, preservation of the geometry under
refinement, exact integration of curved NURBS areas, unit outward boundary
normals and the divergence theorem.

\subsection{Differentiability}
\label{subsec:ex_gradients}

To verify the sensitivities, we compare adjoint gradients with central
finite-difference approximations in Table~\ref{tab:gradients}. The inputs
considered here are material parameters, source terms and the control points
which define the geometry.

\begin{table}[htbp]
\centering
\caption{Adjoint gradients against central finite differences. The finite-difference values include truncation and round-off error.}
\label{tab:gradients}
\small
\begin{tabular}{l l l l}
\hline
differentiated w.r.t. & adjoint & central differences & rel. difference \\
\hline
material $E$ & -7.61302756e-11 & -7.61302756e-11 & 2.2e-10 \\
material $\nu$ & 1.46080230e-06 & 1.46080230e-06 & 2.3e-09 \\
source amplitude & 8.38797967e-02 & 8.38797967e-02 & 8.1e-14 \\
control point & -2.42519379e-02 & -2.42519379e-02 & 1.6e-10 \\
\hline
\end{tabular}
\end{table}

The largest relative difference is $2.3\times10^{-9}$, including the
truncation and round-off errors in the finite-difference approximation. In
Table~\ref{tab:implicit}, we give the corresponding checks for the iterative
solution procedures. By counting the calls to the linear solver, we verify
that each backward pass solves exactly one adjoint system, independently
of the number of parameters or forward iterations. We also observe agreement
to round-off between each gradient and its \texttt{jax.jit} compilation.
For this comparison, the entire solution procedure, including the convergence
tests, must
be represented by a single traced program.

\begin{table}[htbp]
\centering
\caption{Implicit differentiation for each solution method. The nonlinear
Galerkin and energy-minimization runs take five to twelve Newton steps and
several hundred L-BFGS steps, respectively. In each case, the backward pass
requires the solution of one adjoint system.}
\label{tab:implicit}
\small
\begin{tabular}{l l l l c c}
\hline
method & param. & adjoint & central diff. & rel. diff. & adjoint solutions \\
\hline
Galerkin, linear & \texttt{f} & 3.308868e-02 & 3.308868e-02 & 1.9e-11 & 1 \\
Galerkin, Newton & \texttt{load} & 1.340132e-02 & 1.340132e-02 & 3.1e-10 & 1 \\
energy minimization & \texttt{f} & 3.308868e-02 & 3.308868e-02 & 8.6e-11 & 1 \\
energy, nonlinear & \texttt{load} & 1.340132e-02 & 1.340132e-02 & 2.0e-12 & 1 \\
collocation & \texttt{f} & 6.287638e-02 & 6.287638e-02 & 6.8e-12 & 1 \\
collocation, nonlinear & \texttt{g} & 1.236019e+00 & 1.236019e+00 & 4.0e-10 & 1 \\
\hline
\multicolumn{6}{l}{\footnotesize Every gradient also agrees with its \texttt{jax.jit} compilation to a relative difference of} \\
\multicolumn{6}{l}{\footnotesize 2e-13, so the whole solution procedure is one traced program. On the last row the Gauss--Newton} \\
\multicolumn{6}{l}{\footnotesize Jacobian $J^{\top}J$ would give 2.3e-03 relative error instead of 4.0e-10; see Eq.~(\ref{eq:gauss_newton_gap}).} \\
\hline
\end{tabular}
\end{table}

Geometric sensitivities use the same differentiation interface. The
control-point coordinates are stored as traced arrays in the function space,
so the gradient of a compliance functional with respect to them can be
computed directly:

\begin{codebox}\begin{Verbatim}[commandchars=\\\{\}]
\PY{k}{def}\PY{+w}{ }\PY{n+nf}{compliance}\PY{p}{(}\PY{n}{cpts}\PY{p}{)}\PY{p}{:}
    \PY{n}{problem} \PY{o}{=} \PY{n}{base}\PY{o}{.}\PY{n}{with\PYZus{}space}\PY{p}{(}\PY{n}{dataclasses}\PY{o}{.}\PY{n}{replace}\PY{p}{(}\PY{n}{V}\PY{p}{,} \PY{n}{cpts}\PY{o}{=}\PY{n}{cpts}\PY{p}{)}\PY{p}{)}
    \PY{k}{return} \PY{n}{jnp}\PY{o}{.}\PY{n}{sum}\PY{p}{(}\PY{n}{jx}\PY{o}{.}\PY{n}{solve}\PY{p}{(}\PY{n}{problem}\PY{p}{,} \PY{n}{params}\PY{o}{=}\PY{n}{params}\PY{p}{)}\PY{o}{.}\PY{n}{u} \PY{o}{*}\PY{o}{*} \PY{l+m+mi}{2}\PY{p}{)}

\PY{n}{dJ\PYZus{}dcpts} \PY{o}{=} \PY{n}{jax}\PY{o}{.}\PY{n}{grad}\PY{p}{(}\PY{n}{compliance}\PY{p}{)}\PY{p}{(}\PY{n}{V}\PY{o}{.}\PY{n}{cpts}\PY{p}{)}
\end{Verbatim}
\end{codebox}

The same procedure applies to a nonlinear problem, since the implicit
correction~\eqref{eq:implicit} is a traced expression. In this case, the
control points enter the residual through the quadrature data rather than
the parameter argument, but their contribution is still included by automatic
differentiation. We verify the resulting shape derivative for a hyperelastic
problem by comparing it with central finite differences.

\subsection{Computational cost}
\label{subsec:ex_timings}

Table~\ref{tab:timings} separates the cost of setup, solution and gradient
evaluation. The steady-state timings are medians of three repetitions after
warm-up, measured in double precision under the single-core CPU protocol
described at the beginning of Section~\ref{sec:examples}. The cost of the
first call, including JAX tracing and XLA compilation, is reported separately.
When a problem is solved only once, compilation dominates the total cost.

\begin{table}[htbp]
\centering
\caption{Cost of setup, solution and gradient for the two-dimensional Poisson
problem. The ratio of the gradient time to the solution time is close to two,
consistent with one forward and one adjoint solution. The first call includes
compilation.}
\label{tab:timings}
\small
\begin{tabular}{c c r r l l l c c}
\hline
$p$ & ref. & dofs & elems & setup [s] & solution [s] & gradient [s] & $\frac{\text{gradient}}{\text{solution}}$ & first call [s] \\
\hline
2 & 3 & 100 & 64 & 0.017 & 0.182 & 0.343 & 1.88 & 21.23 \\
2 & 4 & 324 & 256 & 0.017 & 0.185 & 0.338 & 1.83 & 21.04 \\
2 & 5 & 1,156 & 1,024 & 0.018 & 0.185 & 0.362 & 1.96 & 21.05 \\
3 & 3 & 121 & 64 & 0.017 & 0.179 & 0.339 & 1.89 & 21.61 \\
3 & 4 & 361 & 256 & 0.017 & 0.183 & 0.346 & 1.89 & 21.80 \\
3 & 5 & 1,225 & 1,024 & 0.018 & 0.189 & 0.363 & 1.91 & 22.33 \\
4 & 3 & 144 & 64 & 0.017 & 0.179 & 0.337 & 1.89 & 21.98 \\
4 & 4 & 400 & 256 & 0.017 & 0.181 & 0.348 & 1.92 & 22.20 \\
4 & 5 & 1,296 & 1,024 & 0.018 & 0.205 & 0.383 & 1.87 & 22.27 \\
\hline
\end{tabular}
\end{table}

First-call compilation dominates at these sizes, and its magnitude varies with
the environment and execution path. The measurements use repeated calls to the
public \texttt{solve} function. Compiling the entire parameter-to-solution
function with \texttt{jax.jit} removes additional Python dispatch and setup
work, so this execution path should be timed separately when assessing a
workload in which the same problem is solved repeatedly.

\subsection{Inverse problem: permeability identification}
\label{subsec:ex_darcy}

As an example of an inverse problem, we consider the recovery of a spatially
varying permeability field from pressure measurements in Darcy flow. The
energy density is
$\psi = \tfrac{1}{2} k(\bm{x}) |\nabla p|^2$ with $k = \exp(\sum_m \theta_m
\phi_m)$ for a truncated cosine basis $\phi_m$, giving nine parameters. Pressures
are sampled at forty interior points from a reference field, and the parameters
are recovered by minimizing the mean squared misfit with Tikhonov
regularization.

\begin{figure}[htbp]
\centering
\includegraphics[width=\textwidth]{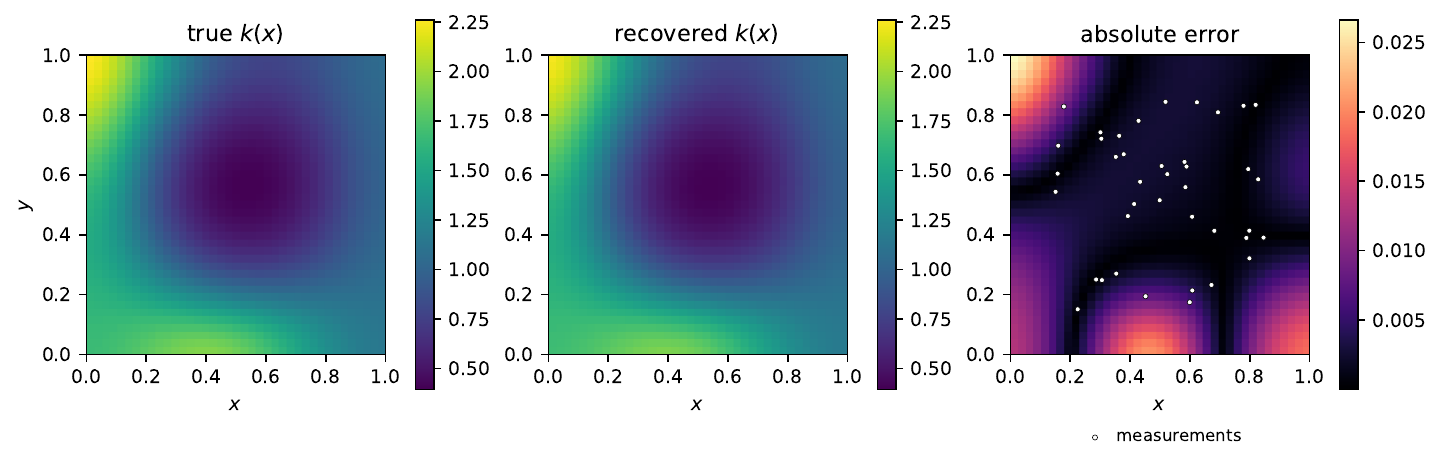}
\caption{Permeability identification. The reference and recovered fields are
shown on the left and in the center, respectively. The right panel shows the
absolute error together with the measurement locations.}
\label{fig:darcy}
\end{figure}

Using L-BFGS, the objective decreases from $7.3\times10^{-5}$ to
$4.1\times10^{-7}$ in thirty iterations. The recovered permeability has a
relative error of $0.59\%$ in the $L^2$ norm, as shown in
Figure~\ref{fig:darcy}. Here the objective includes the solution of the Darcy
problem and the evaluation of pressure at scattered physical points. The
entire computation is compiled as one function and differentiated with
\texttt{jax.grad}, requiring the solution of one adjoint problem for each
gradient evaluation.

\subsection{Material identification from full-field measurements}
\label{subsec:ex_mnist}

We next apply the same procedure to experimental data from the Mechanical MNIST
benchmark~\cite{lejeune2020,mnistchallenge}. Ninety $40\times40\times2$\,mm
specimens were 3D-printed with an MNIST-digit inclusion of one material in a
matrix of another, and stretched in uniaxial tension. Digital image correlation
(DIC) recorded the displacement field on a grid of up to $51\times51$ points.
Four material pairs occur, with inclusion-to-matrix stiffness ratios between
$0.9$ and $11$. The task is to recover the material map from the measured
displacements.

The measured motion of the specimen edges differs from that reported by the
testing machine. In particular, the nominally fixed bottom edge moves by up
to $0.37$\,mm, while the displacement of the top edge is smaller than the
crosshead displacement by up to $0.2$\,mm. We therefore prescribe the measured
edge displacements as Dirichlet data on all four sides. With this choice,
the fixture does not have to be included in the model. Moreover, for linear
elasticity, the interior response depends only on stiffness ratios, so the
unknown absolute modulus and pixel size do not enter the inversion. Since
the Dirichlet values are traced parameters, the same problem definition can
be used for each specimen by changing the data.

We use a single quadratic patch with knots on every interior pixel line,
giving 2,500 elements and 5,408 degrees of freedom. A nodal field $z$ on the
pixel grid defines $\phi=\sigma(z)$, which is interpolated to the quadrature
points using the per-point data described in
Section~\ref{subsec:design_problem}. The stiffness is
$E=E_m\exp(\phi\log r)$ in linear plane stress. We minimize the interior
displacement misfit, weighted per component, plus a total-variation penalty on
$\phi$, using L-BFGS-B for 80 iterations, once for each candidate stiffness ratio
$r$. Each objective and gradient evaluation for the 2,601 parameters requires
one forward and one adjoint solution (Section~\ref{subsec:design_solvers}) and
takes about $0.1$\,s. The gradient agrees with central differences to seven significant digits.

\begin{figure}[htbp]
\centering
\includegraphics[width=\textwidth]{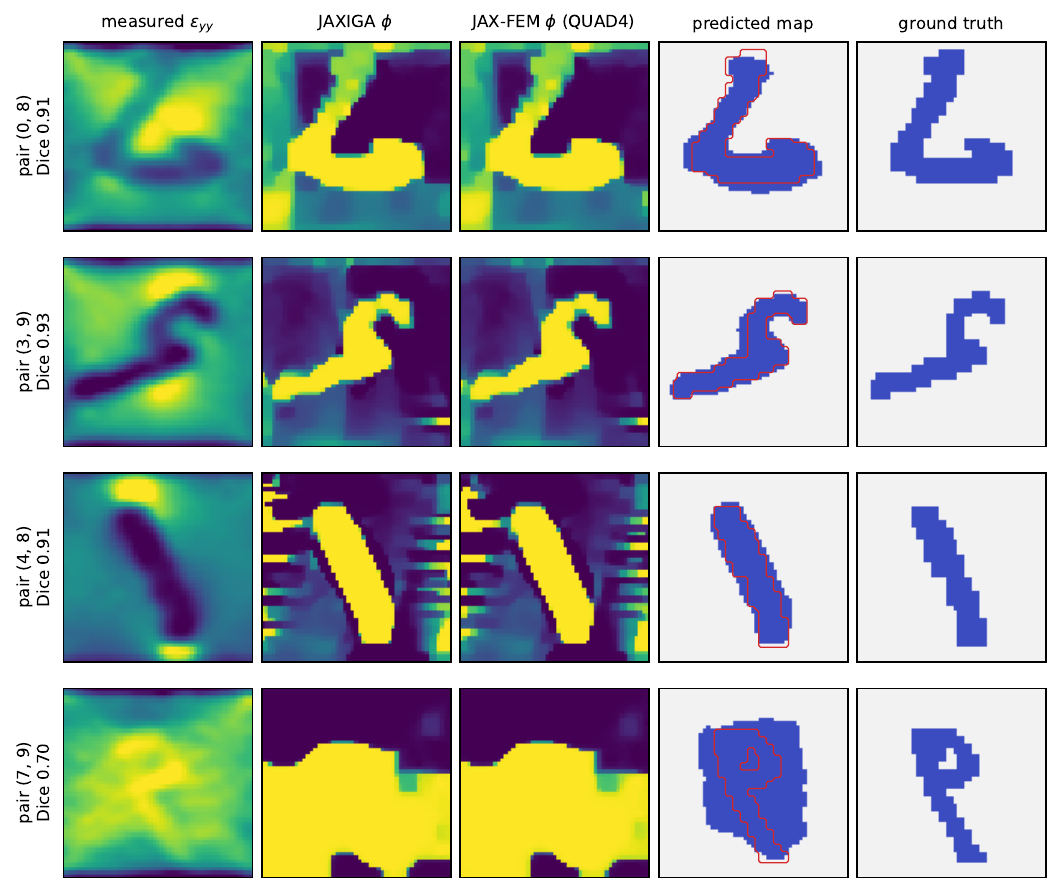}
\caption{Material identification on experimental Mechanical MNIST specimens
using held-out predictions for one median specimen per material pair. From
left to right, the panels show the measured axial strain, the inclusion field
$\phi$ recovered by JAXIGA, the same inversion with JAX-FEM QUAD4 elements,
the predicted material map with the true inclusion outlined in red, and the
reference material map.}
\label{fig:mnist}
\end{figure}

If the recovered field $\phi$ is thresholded directly, the stiff inclusions
are located reliably, but their width is overestimated. This can be explained
by the smoothing in DIC, which limits the available spatial resolution. A
wider inclusion with a more gradual transition at its boundary can give an
equally good fit to the smooth displacement field. We therefore use the four
recovered fields and their final misfits as features in a small learned model,
where they are combined with local strain features for per-pixel logistic
segmentation and material-pair classification.
Under nested five-fold cross-validation over the 90 specimens, adding the
features obtained from the inversions raises the mean Dice coefficient from
$0.801$ to $0.839$ (paired gain
$0.038$, 95\% bootstrap interval $[0.012,0.070]$). Material-pair accuracy rises
from $94.4\%$ to $97.8\%$. Figure~\ref{fig:mnist} shows held-out predictions for
one specimen of each pair. For pair $(7,9)$ the inclusion is only $10\%$ more
compliant than the matrix, and the displacement barely depends on $\phi$.

We repeat the inversion with
JAX-FEM~\cite{jaxfem} QUAD4 elements, with nodes on the DIC points and a similar number of unknowns
(5,202 versus 5,408), and with QUAD8 elements, which match the polynomial degree.
The data, objective, optimizer, quadrature ($3\times3$), SciPy SuperLU solver,
and segmentation and material-pair classification model are identical.
Table~\ref{tab:mnist} and Figure~\ref{fig:mnist} show similar predictive
performance for these discretizations: the Dice coefficient of the final
predicted material map differs by less than $10^{-3}$, no specimen's pair prediction
changes, and the thresholded fields agree on 98--99\% of pixels.
For these computations, JAXIGA has the lowest solution time. On one core of
an Apple M1 Pro CPU, it requires $31$\,s per specimen (four inversions), compared
with $45$\,s for QUAD4 and $132$\,s for QUAD8. On an NVIDIA A10 GPU, the time is
reduced to $16$\,s when Jacobi-preconditioned conjugate gradients are used
and the linear systems are solved on the device. As in
Section~\ref{subsec:ex_parity}, the GPU computations with JAX-FEM are slower
than the CPU computations, both with a direct solver on the host and with its
JAX BiCGSTAB solver. The latter takes
$256$\,s for QUAD4 and does not converge for QUAD8.

The quadratic spline space used here has non-repeated interior knots and is
therefore $C^1$, so it cannot represent a strain jump across a bonded material
interface. For an interface aligned with knot lines, such a jump could be
admitted by increasing the corresponding knot multiplicity to two, reducing
continuity across the interface to $C^0$. To examine whether this continuity
affects the identification, we also solved with a pixel-constant material on a
QUAD4 mesh refined so that every pixel boundary is an element edge. The $C^0$
space then permits a strain jump at the interface. This changes the forward solution by 2--3.5\% of its
geometry-dependent part, but it fits the measured field no better (mean
misfit $0.343$ in both cases). It also lowers the Dice of the thresholded field
from $0.721$ to $0.710$. At this DIC resolution, the comparison does not show an advantage for the
sharp-interface model.

\begin{table}[htbp]
\centering
\caption{Material identification with JAXIGA and JAX-FEM, using the same
procedure except for the discretization. Dice and material-pair accuracy are
reported for the full model under nested cross-validation. ``$\phi$ alone''
denotes the Dice coefficient of the thresholded inclusion field for the true
material pair. Times are given in seconds per specimen (four inversions).
The CPU computations use one core of an Apple M1 Pro CPU with SciPy SuperLU.
The GPU times correspond to the fastest converged solver on an NVIDIA A10 GPU,
which is conjugate gradients for JAXIGA and a direct solver on the host for JAX-FEM.}
\label{tab:mnist}
\small
\setlength{\tabcolsep}{5pt}
\begin{tabular}{l r c c c c r r}
\hline
Discretization & dofs & continuity & Dice & pair acc. & $\phi$ alone & CPU (s) & GPU (s) \\
\hline
JAXIGA, B-spline $p=2$ & 5,408 & $C^1$ & 0.839 & 97.8\% & 0.720 & 30.7 & 16.4 \\
JAX-FEM QUAD4 & 5,202 & $C^0$ & 0.839 & 97.8\% & 0.725 & 45.2 & 123.4 \\
JAX-FEM QUAD8 & 15,402 & $C^0$ & 0.839 & 97.8\% & 0.718 & 132.1 & 222.0 \\
\hline
\end{tabular}

\end{table}

\subsection{Adaptive local refinement on a corner singularity}
\label{subsec:ex_adaptive}

To examine the effect of local refinement on a singular solution, we consider
the L-shaped domain $\Omega = (-1,1)^2 \setminus \left([0,1]\times[-1,0]\right)$
and the harmonic function
\begin{equation}
u(r,\theta) = r^{2/3} \sin\!\left(\tfrac{2}{3}\theta\right),
\qquad \theta \in [0, \tfrac{3\pi}{2}],
\end{equation}
whose gradient is unbounded at the re-entrant corner. Since the regularity is
limited to $u \in H^{1+2/3-\epsilon}$, uniform refinement gives convergence
of order $O(N^{-1/3})$ in the energy norm, independently of the polynomial
degree. Thus, increasing the degree alone cannot improve the asymptotic rate
for this problem. By grading the mesh toward the corner, the optimal rate
$O(N^{-p/2})$ can be recovered. The domain is
built from three unit squares, so three patches meet at the corner and
refinement must be mirrored across two interfaces. At each refinement cycle,
the algorithm solves the discrete problem, evaluates the residual indicator,
and selects elements using D\"orfler bulk marking~\citep{dorfler1996} at
$50\%$. It then refines the marked region.

\begin{figure}[htbp]
\centering
\includegraphics[width=\textwidth]{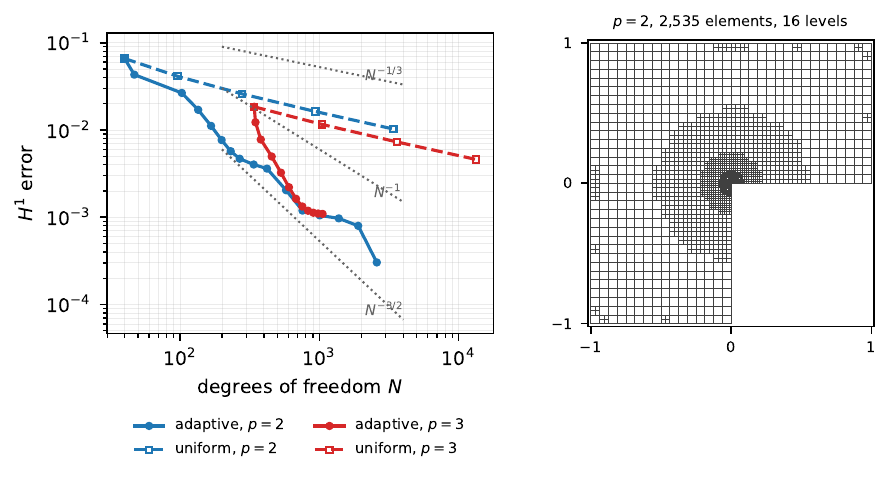}
\caption{Adaptive versus uniform refinement on the L-shaped domain. The left
panel shows energy-norm error against degrees of freedom, with dotted lines
indicating the optimal rates $N^{-p/2}$ and the regularity-limited rate
$N^{-1/3}$. Uniform refinement follows $N^{-1/3}$ for both degrees, so
increasing the degree does not improve the asymptotic rate. Adaptive refinement
recovers the optimal rate. The right panel shows the adapted $p=2$ mesh as
element boxes. The grading into the re-entrant corner spans sixteen levels,
and the three-patch structure is visible in the
interfaces along the negative $x$- and positive $y$-axes.}
\label{fig:lshape}
\end{figure}

\begin{table}[htbp]
\centering
\caption{Adaptive refinement on the L-shaped domain. Rates are least-squares
fits over the final eight cycles. D\"orfler marking refines in bursts, causing
the error to fall in steps. A shorter fitting window would therefore be
sensitive to its position within a step. $n_{\mathrm{loc}}^{\max}$ is the largest number of basis
functions on any element, with the tensor-product value in parentheses.}
\label{tab:lshape}
\begin{tabular}{c c r l r l c c}
\hline
& \multicolumn{3}{c}{adaptive} & \multicolumn{2}{c}{uniform} & & \\
\cmidrule(lr){2-4}\cmidrule(lr){5-6}
$p$ & dofs & $H^1$ error & rate & dofs & $H^1$ error & optimal & $n_{\mathrm{loc}}^{\max}$ \\
\hline
2 & 2,585 & 3.039e-04 & $-1.133$ & 3,400 & 1.025e-02 & $-1.0$ & 12 (9) \\
3 & 1,055 & 1.093e-03 & $-1.556$ & 13,333 & 4.595e-03 & $-1.5$ & 21 (16) \\
\hline
\multicolumn{8}{l}{\footnotesize Uniform refinement rates: $-0.371$ at $p=2$, $-0.363$ at $p=3$, against the regularity limit $-1/3$.} \\
\hline
\end{tabular}
\end{table}

For adaptive refinement, the fitted slopes are $-1.13$ and $-1.56$ for
$p=2,3$, which are close to the optimal values of $-1$ and $-1.5$. For
uniform refinement, the slopes remain near $-1/3$, as shown in
Table~\ref{tab:lshape}. At $p=3$, the adaptive solution is four times more
accurate while using approximately twelve times fewer unknowns. With class-2
admissibility, the maximum number of locally supported functions is 12 and
21 after sixteen and twelve refinement levels, respectively.

We note that the residual estimator does not account for error in the Dirichlet data, which
are interpolated at boundary control points. Here the exact solution is nonzero
on the outer boundary, so the interpolation error places a lower limit on
accuracy that the marking strategy cannot address. At $p=2$ this limit is
low enough not to affect the observed convergence, but at
$p=3$, starting from the same coarse mesh, the error stops decreasing at
$3.7\times10^{-3}$ while the estimator continues to fall. The $p=3$ run in
Table~\ref{tab:lshape} therefore starts from a mesh with two additional uniform
refinements, which resolves the boundary data. Projecting the boundary data and accounting for
their error in the marking strategy are priorities for further development
(Section~\ref{subsec:disc_limitations}).

\subsection{A fourth-order problem: the Kirchhoff--Love shell}
\label{subsec:ex_shell}

As an example of a fourth-order problem, we consider shell bending, whose
primal variational formulation requires $H^2$ conformity because it depends
on second derivatives of the displacement. For a curved midsurface, membrane
and bending action are also coupled through the geometry. We use the
isogeometric Kirchhoff--Love formulation of Kiendl et al.~\cite{kiendl2009},
where the midsurface is represented by a single NURBS patch with control
points in $\mathbb{R}^3$. The displacement $\bm{u}$ is approximated in the
same rational space (\texttt{vec}$\,=3$). The strains for the exact deformed
surface are given by
\begin{equation}
\varepsilon_{\alpha\beta} = \tfrac{1}{2}\left(\bar{\bm{a}}_\alpha\cdot\bar{\bm{a}}_\beta
 - \bm{a}_\alpha\cdot\bm{a}_\beta\right), \qquad
\kappa_{\alpha\beta} = \bm{a}_{\alpha,\beta}\cdot\bm{a}_3
 - \bar{\bm{a}}_{\alpha,\beta}\cdot\bar{\bm{a}}_3,
\end{equation}
where $\bar{\bm{a}}_\alpha = \bm{a}_\alpha + \bm{u}_{,\alpha}$ and $\bar{\bm{a}}_3$
is the unit normal of the deformed surface. After transformation to a local
Cartesian frame, these strains are used in the energy density
$\psi = \tfrac{t}{2}\,\bm{\varepsilon}^{\mathsf T}\mathbf{D}\bm{\varepsilon}
+ \tfrac{t^3}{24}\,\bm{\kappa}^{\mathsf T}\mathbf{D}\bm{\kappa}$, where $\mathbf{D}$ is
the plane-stress elasticity matrix and $t$ is the thickness.

The linear shell uses the exact linearization of these strains at
$\bm{u}=\bm{0}$, obtained by forward-mode differentiation (\texttt{jax.jvp})
rather than from a hand-derived strain--displacement matrix, while the
geometrically nonlinear shell uses the strains themselves. The surface-patch
evaluation is described in Section~\ref{subsec:design_hessian}. Before
considering the benchmark, we perform two consistency checks. On a flat patch,
the shell solution agrees with the plate solution to a relative error of
$4.7\times10^{-15}$, and the clamped square plate
under uniform load gives a center deflection that is $1.00001$ times the
tabulated value $0.0012653\,qa^4/D$~\citep{timoshenko1959}. On the curved roof, an
infinitesimal rigid rotation produces zero linearized strain and a finite
rotation produces zero exact strain.

\begin{figure}[htbp]
\centering
\includegraphics[width=\textwidth]{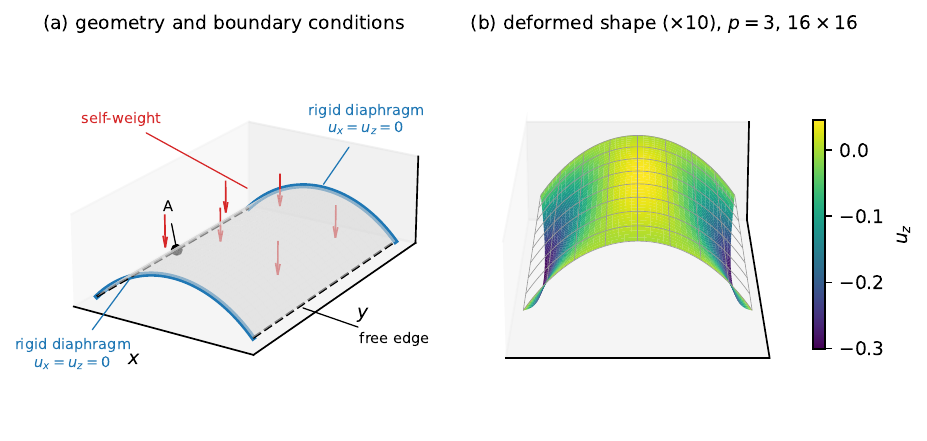}
\caption{Scordelis--Lo roof. (a) Geometry and boundary conditions: cylindrical
shell with radius $R=25$, length $L=50$, opening angle $80^\circ$, thickness
$t=0.25$, $E=4.32\times10^8$ and $\nu=0$, under self-weight 90 per unit area.
Rigid diaphragms at the curved ends impose $u_x=u_z=0$, while the straight edges are
free. Point A is the midpoint of a free edge. (b) Deformed shape viewed along
the roof axis, displacements magnified ten times and colored by the vertical
displacement, with the undeformed midsurface as a gray wireframe. Both free
edges move down and inward symmetrically.}
\label{fig:shell_roof}
\end{figure}

The Scordelis--Lo roof \citep{scordelis1964} is part of the shell obstacle
course \citep{macneal1985,belytschko1985}. It combines membrane and bending
action on a singly curved shell. Figure~\ref{fig:shell_roof}(a) shows the
geometry and supports. The rigid diaphragms leave the axial rigid translation
free, so a single axial displacement component is fixed at one corner to
remove this mode without changing the strains. The midsurface is one quadratic-by-linear NURBS
patch, so the cylinder is exact at every resolution. It is degree-elevated to
$p$ and refined by knot insertion.

\begin{figure}[htbp]
\centering
\includegraphics[width=0.55\textwidth]{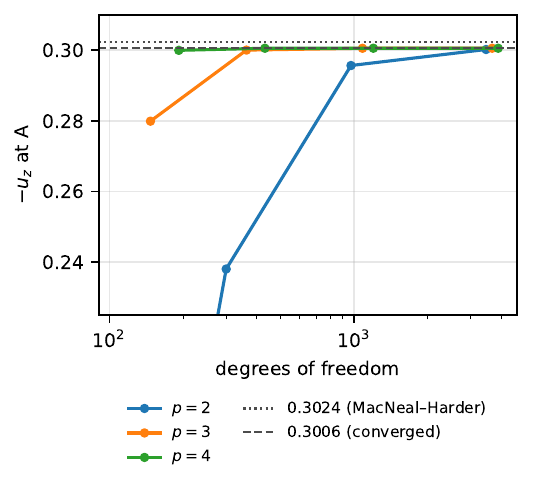}
\caption{Scordelis--Lo roof: convergence of the downward displacement magnitude at A with
the number of degrees of freedom. The $p=2$ value on $4\times4$ elements,
$0.069$, lies below the plotted range.}
\label{fig:shell_convergence}
\end{figure}

The downward displacement magnitude $-u_z$ at A is given in
Table~\ref{tab:shell} and Figure~\ref{fig:shell_convergence}. For $p=4$, we
obtain the converged value $0.30059$ on $8\times8$ elements (432 degrees of
freedom), while $p=3$ requires $16\times16$ elements. The convergence for
$p=2$ is slower, with a displacement of $0.3002$ on $32\times32$ elements.
This behavior is consistent with membrane locking in the quadratic
discretization~\cite{casquero2023cas}. The
converged value agrees with Kirchhoff--Love results for this benchmark in the
isogeometric literature~\citep{kiendl2009}. The frequently quoted reference
$0.3024$~\citep{macneal1985} lies $0.6\%$ higher. Values in this range are
reported for different shell theories.
Figure~\ref{fig:shell_roof}(b) shows the characteristic deformation. At A the
displacement is $(u_x,u_z)=(0.158,-0.301)$, and the opposite edge moves by
$(-0.158,-0.301)$, so the symmetry is preserved to round-off. Both free edges
move downward and inward, while the crown rises slightly.

\begin{table}[htbp]
\centering
\caption{Scordelis--Lo roof: downward displacement magnitude $-u_z$ at the free-edge midpoint A
for degree $p$ and $n\times n$ elements. The converged value is $0.30059$,
compared with the reference value of $0.3024$ \citep{macneal1985}.}
\label{tab:shell}
\small
\begin{tabular}{c r r r r}
\hline
 & \multicolumn{4}{c}{$-u_z$ at A for $n\times n$ elements} \\
$p$ & $4\times4$ & $8\times8$ & $16\times16$ & $32\times32$ \\
\hline
2 & 0.0695 & 0.2381 & 0.2957 & 0.3002 \\
3 & 0.2799 & 0.3001 & 0.3006 & 0.3006 \\
4 & 0.3000 & 0.3006 & 0.3006 & 0.3006 \\
\hline
\end{tabular}
\end{table}

\subsection{Assembly variants, batching and time integration}
\label{subsec:ex_variants}

In this section, we consider two types of problems to examine the capabilities
described in Section~\ref{subsec:design_variants}. We first use a Poisson
problem to compare the assembly and batching procedures. We then consider
longitudinal wave propagation in an elastic bar, both as a forward problem
for checking time integration and as an inverse problem for identifying a
spatially varying stiffness.

\subsubsection*{Assembly and batching checks}

For the first four rows of Table~\ref{tab:capabilities}, we use the
two-dimensional Poisson problem
\begin{equation}
 -\nabla\cdot(a_0\nabla u)=f \quad\text{in }\Omega=(0,1)^2,
 \qquad u=0 \quad\text{on }\partial\Omega,
 \qquad f=8\pi^2\sin(2\pi x)\sin(2\pi y).
 \label{eq:variants_poisson}
\end{equation}
Here $a_0$ is a positive constant, so the exact solution is
$u=\sin(2\pi x)\sin(2\pi y)/a_0$. We use quadratic splines on a uniform
$16\times16$ mesh with $a_0=1$ to compare matrix-free tangent products and
solutions with their assembled counterparts. The same problem is used for
chunked assembly, with eight elements per chunk. For the batching check, we
use an $8\times8$ mesh and eight equally spaced values of $a_0$ between
$0.5$ and $4.0$, while keeping the source term fixed. The solutions computed
using \texttt{vmap} are compared with those obtained by solving each problem
separately. Since these comparisons concern different implementations of the
same discrete problems, agreement to round-off is expected.

\begin{table}[htbp]
\centering
\caption{Consistency of the alternative assembly and solution paths. The first four rows compare alternative implementations, while the last two
check analytical frequencies and discrete energy conservation.}
\label{tab:capabilities}
\begin{tabular}{l l l}
\hline
Check & Discrepancy & Measure \\
\hline
matrix-free matvec vs assembled $Kv$ & 3.7e-16 & relative \\
matrix-free vs assembled solution & 1.1e-14 & absolute \\
chunked vs full assembly (8 elements/chunk) & 8.9e-16 & absolute \\
batched (\texttt{vmap}) vs looped solutions, 8 parameters & 1.3e-15 & absolute \\
bar natural frequencies, first five modes & 2.6e-07 & relative \\
Newmark energy drift, 1,000 steps ($\beta=\tfrac14$, $\gamma=\tfrac12$) & 5.1e-15 & relative \\
\hline
\end{tabular}

\end{table}

\subsubsection*{Free vibration of an elastic bar}

For the dynamic problems, let $u(x,t)$ denote the longitudinal displacement
of an undamped elastic bar. With no body force, the governing equation is
\begin{equation}
 \rho\,\frac{\partial^2 u}{\partial t^2}
 -\frac{\partial}{\partial x}
   \left(E(x)\frac{\partial u}{\partial x}\right)=0,
 \qquad 0<x<L,\quad t>0,
 \label{eq:bar_wave}
\end{equation}
where $\rho$ is the density and $E(x)$ is the elastic modulus. The strain
energy density is $\psi=\tfrac12 E(x)u_{,x}^2$. After the Dirichlet
constraints have been imposed, the Galerkin discretization gives
\begin{equation}
 M\ddot{\bm d}(t)+K\bm d(t)=\bm b(t),\qquad
 M_{ij}=\int_0^L\rho R_iR_j\,dx,\qquad
 K_{ij}=\int_0^L E(x)R_{i,x}R_{j,x}\,dx,
 \label{eq:bar_semidiscrete}
\end{equation}
where $\bm d$ contains the free displacement coefficients, $M$ is the
consistent mass matrix and $\bm b$ contains the applied boundary loads.

For the frequency and energy checks, we take $L=1$, $E=2$ and $\rho=3$ in
consistent model units, with both ends fixed, i.e. $u(0,t)=u(L,t)=0$.
The natural frequencies are obtained from the generalized eigenvalue problem
$K\bm q_k=\omega_k^2M\bm q_k$ and compared with
$\omega_k = k\pi c / L$, where $c=\sqrt{E/\rho}$ and $k=1,\ldots,5$.
The maximum relative error is $2.6\times10^{-7}$ for cubic splines
($p=3$) on 32 uniform elements.

For the energy check, the bar undergoes free vibration, so $\bm b(t)=\bm0$.
The initial free coefficients are set to $d_i(0)=\sin(\pi X_i/L)$, where
$X_i$ is the coordinate of the corresponding control point, and the initial
velocity is zero. Thus, the initial displacement is a spline approximation
of the first sinusoidal mode. We use the Newmark method with $\beta = 1/4$
and $\gamma = 1/2$, a time step of $\Delta t=2\times10^{-3}$ and 1,000 steps.
The discrete energy
$\tfrac12\dot{\bm d}^{\!\top}M\dot{\bm d}+\tfrac12\bm d^{\!\top}K\bm d$
is conserved to a relative difference of $5.1\times10^{-15}$, which is of
the order of round-off error. Since the average-acceleration rule conserves
this energy exactly for a linear undamped system, the result provides a
check of the implementation.

\subsubsection*{Identification of a stiffness inclusion}

We next use~\eqref{eq:bar_wave} to recover a stiffness inclusion from a
single velocity trace. Here $L=1$ and $\rho=1$, and the modulus is
parameterized by
\begin{equation}
 E(x;\bm\theta)=E_0\left[1+\alpha
   \exp\left(-\left(\frac{x-x_c}{w}\right)^2\right)\right],
 \qquad \bm\theta=(\alpha,x_c),\quad E_0=1,\quad w=0.12.
 \label{eq:fwi_modulus}
\end{equation}
The unknown parameters are the contrast $\alpha$ and the center $x_c$ of
the inclusion, while its width $w$ is fixed. The right end is fixed and a
Ricker traction pulse is applied at the left end. The initial and boundary
conditions are
\begin{align}
 u(x,0)&=0, & \dot u(x,0)&=0,\label{eq:fwi_initial}\\
 u(L,t)&=0, & -E(0;\bm\theta)u_{,x}(0,t)&=s(t),\label{eq:fwi_boundary}
\end{align}
where
\begin{equation}
 s(t)=\left[1-2\pi^2 f_0^2(t-t_0)^2\right]
       \exp\left[-\pi^2 f_0^2(t-t_0)^2\right],
 \label{eq:fwi_pulse}
\end{equation}
with $f_0=3.5$ and $t_0=0.35$.
The minus sign in the left boundary condition follows from the outward
normal at $x=0$. We use quadratic splines on 128 uniform elements and the
same Newmark parameters as above, with $\Delta t=5\times10^{-3}$ and 360
steps, giving a final time $T=1.8$. The measured quantity is the velocity at
the loaded end, $v_n(\bm\theta)=\dot u_h(0,t_n)$, with $t_n=n\Delta t$.

Synthetic observations $v_n^{\mathrm{obs}}$ are generated with the same
discretization using $\alpha=0.60$ and $x_c=0.62$, without added noise.
Starting from $\alpha=0.20$ and $x_c=0.50$, we minimize the normalized
least-squares objective
\begin{equation}
 J(\bm\theta)=
 \frac{\displaystyle\sum_{n=0}^{360}
       \left[v_n(\bm\theta)-v_n^{\mathrm{obs}}\right]^2}
      {\displaystyle 2\sum_{n=0}^{360}\left(v_n^{\mathrm{obs}}\right)^2}.
 \label{eq:fwi_objective}
\end{equation}
We use 90 Adam iterations with a learning rate of $0.03$. The gradient is
computed by reverse-mode differentiation through the full time integration
using \texttt{scan}, with checkpointing every 40 steps. A dense linear solver
is used at each time step. Thus, this example differentiates the discrete
trajectory, rather than solving a single static adjoint problem for the final
state.

The difference between the observed trace and the trace for a homogeneous
bar has a norm equal to $7.9\%$ of the observed trace norm. This is the signal
available for identifying the inclusion. The recovered contrast and center
are shown in Table~\ref{tab:fwi}, together with the decrease in the objective.
The last two rows compare the derivatives at the initial guess with central
finite differences using a perturbation of $10^{-5}$ in each parameter.
We note that the misfit
oscillates as the arrival time of the reflected pulse changes. Consequently,
if the initial parameters are too far from the true values, the inversion
can converge to an incorrect cycle of the observed waveform. This phenomenon,
known as cycle skipping, is a main difficulty in full-waveform inversion.
We use a low source frequency in this example, while practical applications
commonly use a sequence of increasing frequencies.

\begin{table}[htbp]
\centering
\caption{One-dimensional full-waveform inversion for the bar
in~\eqref{eq:bar_wave}--\eqref{eq:fwi_pulse}. The contrast and center in
the first two parameter rows are shown as true values followed by recovered
values. The misfit is the objective~\eqref{eq:fwi_objective}, before and
after optimization. Derivatives are evaluated at the initial guess.}
\label{tab:fwi}
\small
\begin{tabular}{l l l}
\hline
Quantity & Value & Note \\
\hline
inclusion echo & 7.9\% of trace norm & the entire signal \\
contrast & 0.600 $\to$ 0.6033 & 0.6\% error \\
center & 0.620 $\to$ 0.6199 & $<0.1\%$ error \\
misfit & 4.569e-03 $\to$ 1.774e-07 & 25,757$\times$ lower \\
\hline
$\partial$misfit$/\partial$center & -2.497087e-02 & vs.\ central diff., 6e-09 rel. \\
$\partial$misfit$/\partial$contrast & 1.055024e-02 & vs.\ central diff., 6e-11 rel. \\
\hline
\end{tabular}

\end{table}

\begin{figure}[htbp]
\centering
\includegraphics[width=\textwidth]{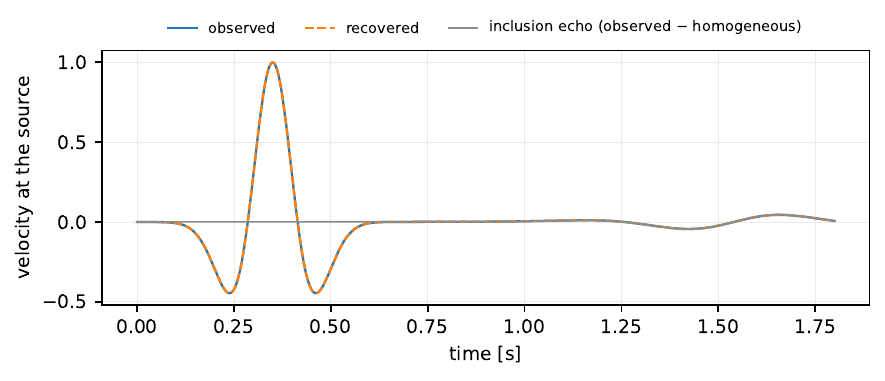}
\caption{Velocity at the loaded end $x=0$ of the elastic bar. The observed and recovered traces are
indistinguishable at this scale. The gray curve shows the difference between
the observed trace and that of a homogeneous bar. This difference isolates the
inclusion's contribution, which supplies the information used in the inversion.}
\label{fig:fwi}
\end{figure}

\subsection{Phase-field fracture}
\label{subsec:ex_phase_field}

Finally, we consider a coupled problem in which the solution depends on the
loading history. The domain is a $1\times1$\,mm plate in plane strain, with
$E=210$\,GPa, $\nu=0.3$, $G_c=2.7$\,N/mm and $\ell=0.0125$\,mm. The bottom
edge is fixed, and a vertical displacement is prescribed on the top edge,
where the horizontal displacement is zero. An initial band in the history
field defines a crack extending from the left edge to the center. Thus, the
notch is represented by damage and is not removed from the geometry.

\subsubsection*{Formulation}

In a phase-field model, the scalar damage variable $\phi:\Omega\to[0,1]$ is
zero for intact material and one for fully broken material. We use an AT2
crack-density term with isotropic stiffness degradation, while damage is
driven by the tensile-energy history. This is a hybrid treatment related to
that of Ambati et al.~\cite{ambati2015hybrid}. The underlying elastic and
fracture energy is given by
\begin{equation}
\mathcal{E}(\mathbf{u}, \phi) = \int_\Omega \left[ g(\phi)\, \psi(\boldsymbol{\varepsilon}) + G_c \left( \frac{\phi^2}{2\ell} + \frac{\ell}{2} |\nabla \phi|^2 \right) \right] d\Omega - W_{\text{ext}}(\mathbf{u}),
\label{eq:phase_field_energy}
\end{equation}
where $g(\phi) = (1-\phi)^2 + k$ is the degradation function, with residual
stiffness $k \ll 1$, and $\psi = \tfrac12 \boldsymbol{\varepsilon} : \mathbb{C} :
\boldsymbol{\varepsilon}$ is the elastic strain energy density. Here $G_c$
denotes the critical energy release rate and $\ell$ is the regularization
length. Irreversibility and growth under tension are introduced through the
driving force of the damage subproblem. In particular, we use the history field
\begin{equation}
\mathcal{H}(\bm{x}, t) = \max_{s \le t} \psi^{+}\!\left(\boldsymbol{\varepsilon}(\bm{x}, s)\right),
\label{eq:history}
\end{equation}
which records the running maximum of the tensile strain energy. Here,
$\psi^{+}$ is obtained from the spectral decomposition of the
strain~\citep{miehe2010}. The maximum-history
construction supplies irreversibility, although pointwise monotonicity after
discretization and mesh transfer must still be checked. Degrading the whole
energy rather than only $\psi^{+}$ keeps the displacement subproblem linear at
frozen $\phi$, since the split then does not enter the tangent. With this
treatment, a fully broken crack carries no compressive stiffness across its faces,
so crack closure is not represented.

With the damage and history fields held fixed in their respective
subproblems, the weak equations are
\begin{align}
 \int_\Omega g(\phi)\mathbb C:\bm\varepsilon(\bm u):
 \bm\varepsilon(\bm v)\,d\Omega &= W_{\mathrm{ext}}(\bm v),\\
 \int_\Omega \left[(G_c/\ell+2\mathcal H)\phi q
             +G_c\ell\nabla\phi\cdot\nabla q\right]d\Omega
 &=\int_\Omega 2\mathcal Hq\,d\Omega .
\end{align}
Both fields use the same spline space. The load loop alternates between
solving these two equations and updating~\eqref{eq:history}. This produces a staggered fixed-point scheme
that does not, in general, correspond to coordinate descent
of~\eqref{eq:phase_field_energy}. In particular,
$\mathcal H=\psi^+$ does not recover that energy's damage equation when
$\psi^-$ is nonzero.

\subsubsection*{Discretization and staggered solution}

We use quadratic THB-splines on a $20\times20$ base mesh, with four levels of
refinement along the crack path. The refined elements remain square and give
$\ell/h=4$. The maximum load depends substantially on the resolution. Over
the prescribed loading range, the maximum loads are $771$, $641$ and
$616$\,N for $\ell/h=1,2,4$, respectively, and the coarsest discretization
does not fracture within this range. We therefore do not regard the finest
curve alone as evidence of spatial convergence.

The load schedule uses 44 increments of $10^{-4}$\,mm followed by 1,156 of
$10^{-6}$\,mm. The staggered driver checks the relative damage residual at
the incoming field against $10^{-4}$ and records convergence flags. The
computation takes 2,422 sweeps over 1,200 increments. The basis data and
sparsity patterns are reused to avoid repeated setup. SuperLU factors are reused only when
iterative refinement meets a relative linear residual tolerance of $10^{-5}$.
Reactions are evaluated from the constrained displacement residual.

\begin{figure}[htbp]
\centering
\includegraphics[width=\textwidth]{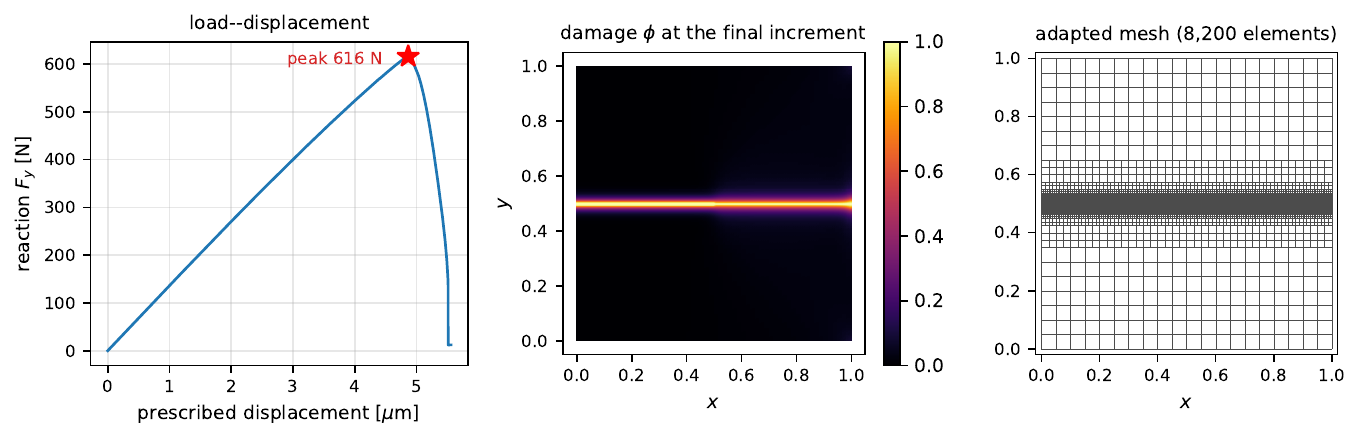}
\caption{Single-edge-notched plate in tension. The left panel shows the
load--displacement response over $1{,}200$ fixed increments. The reaction rises
elastically to $616$~N and then falls to about $12$~N as the crack traverses
the ligament. The center panel shows the damage field at the final increment.
The right panel shows the mesh, refined hierarchically along the crack path
from a $20\times20$ square base so that elements stay
square at $\ell/h = 4$ in both directions.}
\label{fig:phase_field}
\end{figure}

\subsubsection*{Results}

The load--displacement response and final damage field are shown in
Figure~\ref{fig:phase_field}. We observe that the reaction first increases
approximately linearly, reaching a maximum of $616$~N at a prescribed
displacement of $4.86\times10^{-3}$~mm. It then decreases to about $12$~N at
the final increment. The peak
lies within $0.7\%$ of the $620$~N at $5.07\times10^{-3}$~mm reported by the
reference implementation~\cite{igapack} on an adaptively refined cubic mesh with
the same $G_c$, $\ell$ and elastic constants. However, the reference code
contains both second- and fourth-order phase-field models, each with two
degradation functions, and the distributed load--displacement data do not
identify the model used for this curve. Our formulation is second order
with quadratic degradation and gives a steeper softening branch. Since the
reference model is not known, this difference cannot be attributed to mesh
resolution alone. A quantitative validation would require the same model in
both computations, together with separate studies of load-step and spatial
convergence.

\subsubsection*{Crack-following adaptivity}

When the crack path is not known in advance, we refine elements for which the
damage exceeds $1/2$ at a quadrature point, up to a prescribed maximum level.
A surrounding region of radius $\ell$ in two dimensions or $1.5\ell$ in
three dimensions is refined at the same time. This allows the crack to grow
further before the discretization must be rebuilt. Elements already at the
finest level are excluded from the refinement trigger, and the state is
transferred as described in Section~\ref{subsec:design_hierarchical}.

\begin{figure}[htbp]
\centering
\includegraphics[width=\textwidth]{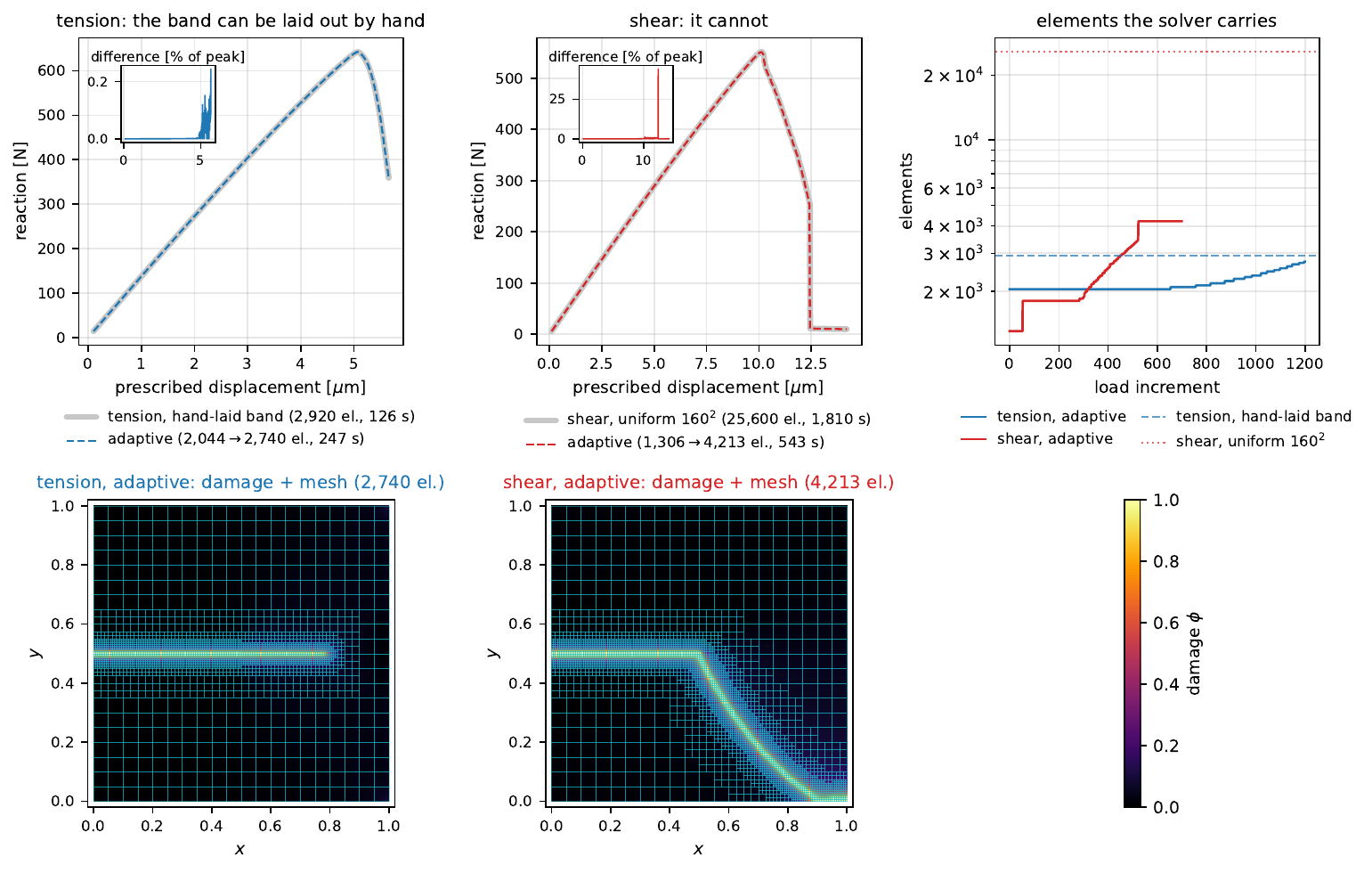}
\caption{Comparison of crack-following adaptivity with fixed meshes. The
upper left panel shows the tension test, for which the crack path is known
and a band of half-width $2\ell$ can be refined in advance. The upper center
panel shows the shear test, where the crack curves toward the lower right
corner. Here a uniform mesh at the same $h$ provides the reference without
requiring prior knowledge of the crack path. The insets show
$100|F_{\mathrm{adaptive}}(u)-F_{\mathrm{reference}}(u)|/F_{\mathrm{reference}}^{\max}$,
comparing reactions at the same prescribed displacement. The shear spike reaches
$44\%$ because the adaptive load drop precedes the uniform one by two increments.
This measures a difference in reactions at a fixed displacement, not a difference
between peak loads. The upper right panel gives the number of elements during
the computation, together with the fixed-mesh values (dashed and dotted).
The lower panels show the final damage field and adapted mesh (cyan) for
tension on the left and shear on the right.}
\label{fig:adaptive_fracture}
\end{figure}

Figure~\ref{fig:adaptive_fracture} compares this strategy with fixed meshes.
The tension computations use $\ell/h=2$, which accounts for the higher peak
load compared with the $\ell/h=4$ case above. For shear, we use
$\ell=0.015$\,mm and $\ell/h=2.4$. The top edge moves horizontally, the
vertical displacement is fixed on the lateral edges, and the bottom edge
remains fixed. The shear loading is applied in
80 increments of $10^{-4}$\,mm followed by 620 of $10^{-5}$\,mm.

For tension, the adaptive and pre-refined computations give peak loads which
differ by $0.23$\,N ($642.74$ versus $642.51$\,N). The adaptive mesh has
2,740 elements at the final increment, compared with 2,920 for the pre-refined
mesh. However, the adaptive computation takes $247$\,s, compared with
$126$\,s for the fixed mesh, because of the cost of rebuilding the
discretization. For shear, the curved crack path requires a uniform reference
mesh with $25{,}600$ elements. The adaptive mesh grows from 1,306 to 4,213
elements, and the peak loads agree within $0.1\%$ ($551.06$ versus
$551.59$\,N). The respective times are $543$ and $1{,}810$\,s, giving a
reduction by a factor of $3.3\times$. 

The shear inset shows a difference which is not apparent in the nearly
vertical part of the load--displacement curve. The adaptive reaction drops at
$u=12.43$\,$\mu$m, while the uniform reaction
drops at $12.45$\,$\mu$m, two increments later. At $12.44$\,$\mu$m the reactions
are $11.25$ and $253.91$\,N, respectively, giving a difference of $43.99\%$ of the
uniform peak, $551.59$\,N. Excluding these two increments, the largest shear
discrepancy is $1.18\%$ of that peak. Thus agreement of the peak loads does not
imply pointwise agreement through abrupt failure. Assessing the shift in the
displacement at which the load drops requires mesh and load-step convergence
studies.

\subsubsection*{Three dimensions}

The same state-transfer procedure and staggered solution method can be used
in three dimensions. For the tensile-energy split, we use regularized
closed-form expressions for the principal strains, which agree with a
symmetric eigensolver to $10^{-13}$. The cost of sparse factorization due to
fill-in can be avoided by using Jacobi-preconditioned CG, as shown in
Table~\ref{tab:cg_vs_direct}. Positive residual stiffness and adequate
constraints ensure positive definiteness. However, these conditions do not
ensure that the conditioning or number of iterations is independent of the
mesh size.

\begin{table}[htbp]
\centering
\caption{Comparison of direct factorization and Jacobi-preconditioned
conjugate gradients for the displacement subproblem of the notched cube.
Fill-in is defined as the ratio of factor nonzeros to matrix nonzeros.
Assembly takes about half a second in these cases, while factorization takes
over twenty times as long at the largest problem size.}
\label{tab:cg_vs_direct}
\small
\begin{tabular}{r r r r r}
\hline
Free dofs & \texttt{splu} & Fill-in & Jacobi-CG & Iterations \\
\hline
$3{,}068$ & $0.25$~s & $2.8\times$ & $0.36$~s & $97$ \\
$7{,}162$ & $1.57$~s & $4.9\times$ & $1.15$~s & $127$ \\
$15{,}452$ & $10.83$~s & $8.5\times$ & $3.84$~s & $157$ \\
\hline
\end{tabular}
\end{table}

\begin{figure}[htbp]
\centering
\includegraphics[width=\textwidth]{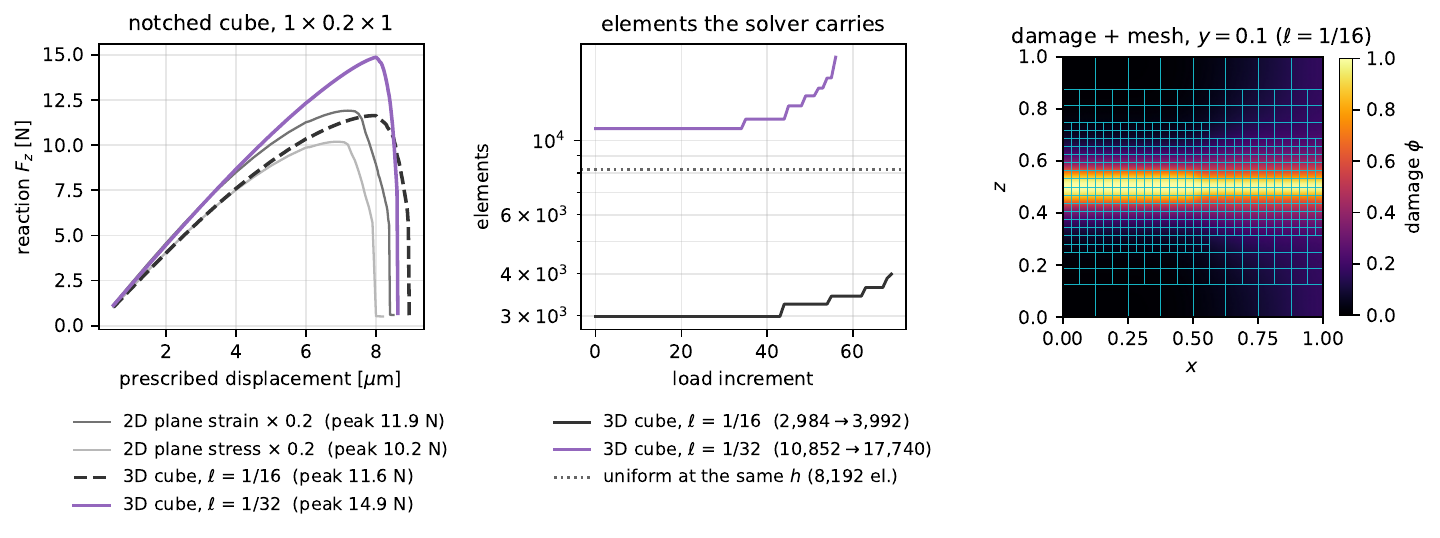}
\caption{Notched cube in tension, with dimensions $1\times0.2\times1$ and a
crack through the thickness. The left panel shows the load--displacement
curves for two regularization lengths. The three-dimensional result lies
between the two-dimensional limits after scaling by the thickness. The
center panel compares the number of elements during adaptive refinement with
a uniform mesh at the same $h$. The right panel shows the final damage field
on the mid-thickness slice $y=0.1$, together with the adapted mesh (cyan),
for $\ell=1/16$.}
\label{fig:cube3d}
\end{figure}

The $1\times0.2\times1$\,mm block has a seeded crack through its thickness
and traction-free lateral faces. At $\ell=1/16$\,mm, its $11.63$\,N peak
lies between the thickness-scaled plane-stress and plane-strain results,
$10.19$ and $11.91$\,N (Figure~\ref{fig:cube3d}). 

When $\ell$ is halved while keeping $\ell/h=2$, the peak increases to
$14.87$\,N, indicating that the seeded band remains dependent on the
regularization length. The finer computation reaches 17,740 elements and
56,112 displacement unknowns over 335 staggered sweeps. Remeshing, setup and
recompilation require $742$ seconds out of the total computation time of
$2{,}083$\,s ($36\%$ of the total). 

\section{Discussion}
\label{sec:discussion}
\label{subsec:disc_limitations}
\label{subsec:disc_future}

The numerical examples illustrate the circumstances in which the present
implementation can reduce computational cost. In particular, when a problem is
solved repeatedly on a fixed mesh, compilation can be reused, while implicit gradients and
linear algebra performed on the device reduce the cost of sensitivity
analysis and solution. Adaptive refinement can reduce the number of unknowns,
but the cost of rebuilding the discretization and recompiling the computation
must also be taken into account. Several improvements in accuracy remain to
be investigated. These include the projection of boundary data, error
indicators which account for boundary and interface residuals, and fracture
convergence studies in which the same model is used by both implementations.

We also mention some limitations of the current implementation. Multipatch
coupling requires conforming interfaces and provides only $C^0$ continuity.
Therefore, primal fourth-order problems must be represented by suitable
single patches. Smooth spline constructions on unstructured meshes, such as
blended B-splines~\cite{wei2018blended} and analysis-suitable unstructured
T-splines~\cite{casquero2020seamless,wei2022asut}, offer a possible way of
removing this restriction for geometries of arbitrary topology. Collocation is
restricted to tensor-product spaces and uses dense normal equations. Its
extension to locally refined spaces could follow the approach developed for
analysis-suitable T-splines~\cite{casquero2016collocation}. Local hierarchical refinement is currently unsupported
for spaces with collapsed-boundary constraints, which tie coefficients together
to ensure that the solution is single-valued on collapsed patch boundaries.
Tensor-product refinement can be performed before constructing these spaces,
with the required coefficient constraints imposed on the refined discretization.
Further reductions in cost may be
obtained with stronger preconditioners or assembly methods which require
less local storage. Such comparisons should be made at the same solution and
gradient accuracy. Weighted quadrature, sum factorization and low-rank
methods provide established approaches for improving
assembly~\cite{calabro2017,bressan2019,lowrank2015,sangalli2018}. For all such
computations, the sensitivities remain subject to the conditions given in
Section~\ref{subsec:design_solvers}.

\section{Conclusions}
\label{sec:conclusions}

In this paper, we have presented JAXIGA and described how three numerical
methods can be obtained from a common energy density, with spline geometry
and first-order implicit sensitivities. We have compared the computed
solutions with analytical results and independent implementations, and have
verified the sensitivities against finite-difference approximations. We have
also shown how hierarchical extraction allows local refinement of the
variational methods, and how smooth single-patch spaces can be used for
fourth-order problems. For the reported plate benchmark, GPU assembly is accelerated
by $24\times$, while the total solution time is reduced by $23\times$
relative to the reported CPU direct solver when a suitable device solver is
used. Thus, the reduction in assembly cost must be accompanied by an
appropriate choice of solver in order to obtain a comparable reduction in
the total computational cost.

\section*{Acknowledgments}

This work was supported in part by the US National Science Foundation (NSF:
\#1953323 and \#2332084), the Advanced Research Projects Agency for Health
(ARPA-H, D25AC00460), U.S. ARMY Medical Research Acquisition Activity (USAMRAA,
HT94252510821), the Manufacturing PA Innovation Program and CMU Digital Twin
for Complex Engineering Systems Initiative.

\section*{Data Availability}

The JAXIGA source code is available under the MIT license at
\url{https://github.com/CMU-CBML/JAXIGA/}.

\appendix

\section{Installation and Usage}
\label{app:installation}

The library requires Python 3.10 or later, together with JAX, NumPy, SciPy
and \texttt{optax}. The dependencies for VTK output, plotting and Gmsh are
optional. After obtaining a repository checkout, the library can be installed
with the following command:

\begin{codebox}\begin{Verbatim}[commandchars=\\\{\}]
pip\PY{+w}{ }install\PY{+w}{ }\PYZhy{}e\PY{+w}{ }\PY{l+s+s2}{\PYZdq{}.[post,plot]\PYZdq{}}
\end{Verbatim}
\end{codebox}

Sixty-four-bit floating point is enabled on import to support the accuracy
needed for the reported discretization and gradient checks. The complete Poisson
example below uses Galerkin. To apply energy minimization or collocation to
the same problem, the user changes the \texttt{method} argument.

\begin{codebox}\begin{Verbatim}[commandchars=\\\{\}]
\PY{k+kn}{import}\PY{+w}{ }\PY{n+nn}{jax}\PY{n+nn}{.}\PY{n+nn}{numpy}\PY{+w}{ }\PY{k}{as}\PY{+w}{ }\PY{n+nn}{jnp}
\PY{k+kn}{import}\PY{+w}{ }\PY{n+nn}{jaxiga}\PY{+w}{ }\PY{k}{as}\PY{+w}{ }\PY{n+nn}{jx}

\PY{n}{patch} \PY{o}{=} \PY{n}{jx}\PY{o}{.}\PY{n}{primitives}\PY{o}{.}\PY{n}{quadrilateral}\PY{p}{(}
    \PY{p}{[}\PY{p}{[}\PY{l+m+mi}{0}\PY{p}{,} \PY{l+m+mi}{0}\PY{p}{]}\PY{p}{,} \PY{p}{[}\PY{l+m+mi}{1}\PY{p}{,} \PY{l+m+mi}{0}\PY{p}{]}\PY{p}{,} \PY{p}{[}\PY{l+m+mi}{1}\PY{p}{,} \PY{l+m+mi}{1}\PY{p}{]}\PY{p}{,} \PY{p}{[}\PY{l+m+mi}{0}\PY{p}{,} \PY{l+m+mi}{1}\PY{p}{]}\PY{p}{]}
\PY{p}{)}\PY{o}{.}\PY{n}{elevate}\PY{p}{(}\PY{l+m+mi}{3}\PY{p}{)}\PY{o}{.}\PY{n}{refine}\PY{p}{(}\PY{l+m+mi}{5}\PY{p}{)}

\PY{n}{V} \PY{o}{=} \PY{n}{jx}\PY{o}{.}\PY{n}{FunctionSpace}\PY{p}{(}\PY{n}{patch}\PY{p}{)}

\PY{k}{def}\PY{+w}{ }\PY{n+nf}{source}\PY{p}{(}\PY{n}{x}\PY{p}{,} \PY{n}{p}\PY{p}{)}\PY{p}{:}
    \PY{k}{return} \PY{p}{(}\PY{l+m+mi}{8} \PY{o}{*} \PY{n}{jnp}\PY{o}{.}\PY{n}{pi}\PY{o}{*}\PY{o}{*}\PY{l+m+mi}{2}
            \PY{o}{*} \PY{n}{jnp}\PY{o}{.}\PY{n}{sin}\PY{p}{(}\PY{l+m+mi}{2} \PY{o}{*} \PY{n}{jnp}\PY{o}{.}\PY{n}{pi} \PY{o}{*} \PY{n}{x}\PY{p}{[}\PY{l+m+mi}{0}\PY{p}{]}\PY{p}{)}
            \PY{o}{*} \PY{n}{jnp}\PY{o}{.}\PY{n}{sin}\PY{p}{(}\PY{l+m+mi}{2} \PY{o}{*} \PY{n}{jnp}\PY{o}{.}\PY{n}{pi} \PY{o}{*} \PY{n}{x}\PY{p}{[}\PY{l+m+mi}{1}\PY{p}{]}\PY{p}{)}\PY{p}{)}

\PY{n}{everywhere} \PY{o}{=} \PY{k}{lambda} \PY{n}{x}\PY{p}{:} \PY{n}{jnp}\PY{o}{.}\PY{n}{full}\PY{p}{(}\PY{n}{x}\PY{o}{.}\PY{n}{shape}\PY{p}{[}\PY{l+m+mi}{0}\PY{p}{]}\PY{p}{,} \PY{k+kc}{True}\PY{p}{)}

\PY{n}{problem} \PY{o}{=} \PY{n}{jx}\PY{o}{.}\PY{n}{Poisson}\PY{p}{(}
    \PY{n}{V}\PY{p}{,}
    \PY{n}{dirichlet}\PY{o}{=}\PY{p}{[}\PY{n}{jx}\PY{o}{.}\PY{n}{DirichletBC}\PY{p}{(}\PY{l+m+mf}{0.0}\PY{p}{,} \PY{n}{where}\PY{o}{=}\PY{n}{everywhere}\PY{p}{)}\PY{p}{]}\PY{p}{,}
    \PY{n}{source}\PY{o}{=}\PY{n}{source}\PY{p}{,}
\PY{p}{)}

\PY{n}{sol} \PY{o}{=} \PY{n}{jx}\PY{o}{.}\PY{n}{solve}\PY{p}{(}\PY{n}{problem}\PY{p}{,} \PY{n}{params}\PY{o}{=}\PY{p}{\PYZob{}}\PY{l+s+s2}{\PYZdq{}}\PY{l+s+s2}{a0}\PY{l+s+s2}{\PYZdq{}}\PY{p}{:} \PY{l+m+mf}{1.0}\PY{p}{\PYZcb{}}\PY{p}{)}
\PY{n+nb}{print}\PY{p}{(}\PY{n}{jx}\PY{o}{.}\PY{n}{errornorm}\PY{p}{(}\PY{n}{sol}\PY{p}{,} \PY{n}{exact}\PY{p}{,} \PY{l+s+s2}{\PYZdq{}}\PY{l+s+s2}{L2}\PY{l+s+s2}{\PYZdq{}}\PY{p}{)}\PY{p}{)}
\PY{n}{sol}\PY{o}{.}\PY{n}{to\PYZus{}vtk}\PY{p}{(}\PY{l+s+s2}{\PYZdq{}}\PY{l+s+s2}{poisson.vtu}\PY{l+s+s2}{\PYZdq{}}\PY{p}{)}
\end{Verbatim}
\end{codebox}

\section{Geometry Import and Refinement}
\label{app:geometry}

For two-dimensional CAD input, quadrilateral meshing is first attempted with
Gmsh's quasi-structured algorithm, followed by Blossom-Quad recombination and
then simple recombination. Any remaining triangles can be subdivided into
quadrilaterals. In three dimensions, structured or extruded hexahedral meshes
are retained. For general input, tetrahedra are subdivided into hexahedra,
without assuming that a general-purpose all-hexahedral mesh generator is
available. The automatic sizing procedure repeatedly halves an initial target
size until the cell-quality and sampled CAD-distance checks are satisfied,
subject to a limit on the number of attempts. We note that this heuristic
does not find the minimum number of cells or provide a certified bound on
the geometric error.

High-order conversion solves $V_B C=V_L X$, where $X$ contains Gmsh nodal
coordinates and $V_L,V_B$ evaluate Lagrange and Bernstein bases at an unisolvent
tensor grid. CAD-aware enrichment instead projects boundary Greville samples
and solves for the control points from these projected samples. Polynomial
cell checks combine sampled distortion with Bernstein bounds on the Jacobian
determinant. Cells are rejected when these bounds are inconclusive, and shared
traces must remain compatible.
The option \texttt{rational\_geometry=True} reconstructs supported conic and
cylindrical maps, which ordinary homogeneous refinement then preserves exactly.

The \texttt{imported\_geometry} module in \texttt{paper.experiments} reproduces
the imported-geometry computations and figure. Standalone scripts
in \texttt{jaxiga.examples} include \texttt{gmsh\_plate},
\texttt{connecting\_rod\_3d} and \texttt{blade\_3d}, with mesh/stress plots and
VTK export. Native NURBS data are bundled, so these examples need neither MATLAB
nor an IGAPack checkout.

\end{document}